\documentclass[11pt]{articleplus}   

\usepackage{amssymb,amstext,amsmath,amsthm,bm}

\usepackage{graphicx,wrapfig}
\usepackage[hang]{subfigure}
\usepackage[all]{xy}

\usepackage{url}

\usepackage[numbers,sort]{natbib}

\usepackage[medium,compact]{titlesec}
\usepackage[margin=1 in,nohead]{geometry}

\usepackage{theomac}

\newcommand{\Bvector}[2]{\stackrel{#1}{\mathstrut #2}}

\begin{document}

\title{Periodic body-and-bar frameworks}

\author{\setcounter{footnote}{0}%
\def\thefootnote{\arabic{footnote}}
Ciprian~Borcea$^{1,4}$ 
\and
Ileana~Streinu$^{2,4}$
\and
Shin-ichi Tanigawa$^3$%
}

\maketitle

{
\def\thefootnote{\arabic{footnote}}
\footnotetext[1]{Department of Mathematics, Rider University, Lawrenceville, NJ 08648. {\it email:}
\url{borcea@rider.edu}}
\footnotetext[2]{Department of Computer Science, Smith College, Northampton, MA 01063, USA.
{\it email:} \url{istreinu@smith.edu, streinu@cs.smith.edu}. {\it URL:} \url{http://cs.smith.edu/~streinu}. } 
\footnotetext[3]{Research Institute for Mathematical Sciences (RIMS), Kyoto University, Kyoto, Japan.
{\it email:} \url{tanigawa@kurims.kyoto-u.ac.jp}. 
} 
\footnotetext[4]{Research
supported by a DARPA "23 Mathematical Challenges" grant, under ``Algorithmic Origami and Biology''. }
\setcounter{footnote}{5}
\def\thefootnote{\arabic{footnote}}

\begin{abstract}
	Abstractions of crystalline materials known as periodic body-and-bar frameworks are made of rigid bodies connected by fixed-length bars and subject to the action of a group of translations.  In this paper, we give a Maxwell-Laman characterization for generic minimally rigid periodic body-and-bar frameworks. As a consequence we obtain efficient polynomial time algorithms for their recognition based on matroid partition and pebble games.
\end{abstract}

\section{Introduction}
\label{sec:intro}

In this paper, we study $d$-dimensional periodic frameworks made of rigid bodies connected with rigid bars. Fig.~\ref{fig:prussianBlue} gives an example.  We prove a combinatorial  characterization for the quotient graphs of {\em generic} minimally rigid frameworks, in terms of matroid unions of  graphs satisfying Maxwell-sparsity conditions. 
As a consequence, we obtain an efficient, polynomial time algorithm for their recognition, based on matroid partition and pebble games.

\begin{figure}[h]
    \centering
    {\includegraphics[width=0.24\textwidth]{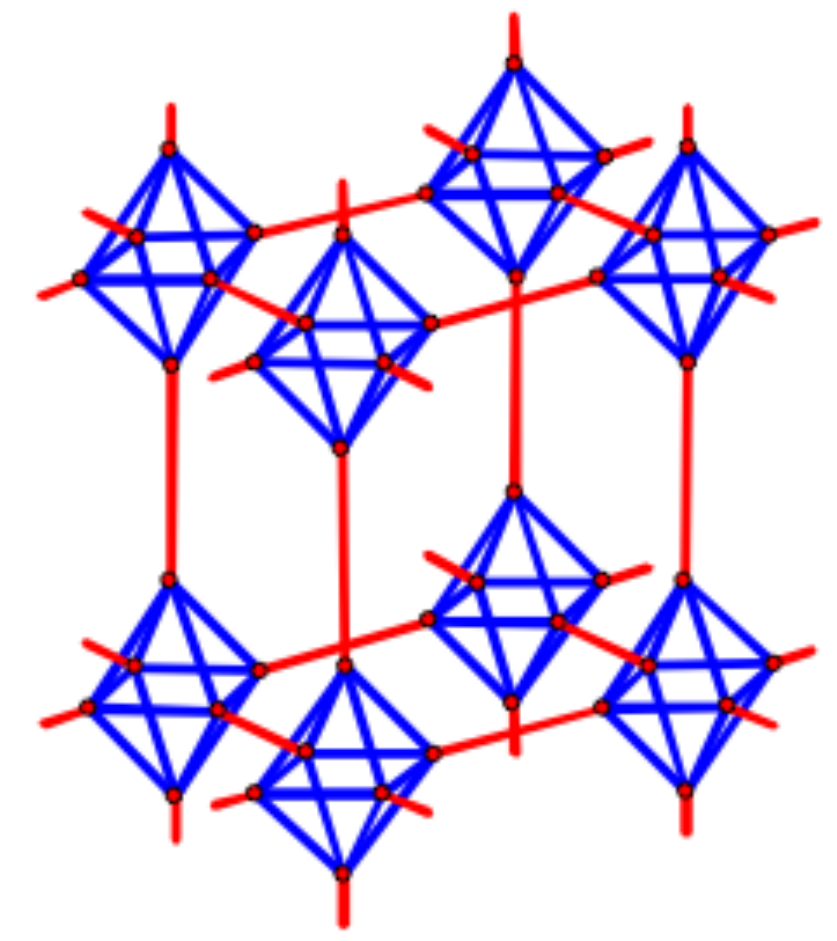}}
    \caption{\small{Fragment of a periodic body-and-bar framework associated to Prussian Blue crystalline materials \cite{goodwin:rums:2006}.
}}
\label{fig:prussianBlue}
\end{figure}

\noindent
{\bf Crystalline materials.} 
{\em Periodic structures} are naturally associated with {\em crystallographic studies} \cite{limaDeFaris:atlasCrystallography:1990}. The central role played by a {\em periodicity lattice} and corresponding fundamental domains was recognized even before X-ray diffraction opened up experimental possibilities. Crystals known as {\em framework materials} \cite{megaw:crystalStructures:1973,dove:displacive:1997,okeeffe:yaghi:etAl:frameworks:2000} have distinctive substructures which can be modeled and investigated as articulated systems made of rigid pieces and joints described abstractly by an {\em infinite periodic graph}. Similar models appear in structural engineering in connection with infinite trusses, foams or cellular materials \cite{deshpande:topologyBending:2001,guest:hutchinson:repetitiveStructures:2003,donev:torquato:energyEfficient:2003,hutchinson:fleck:structuralTruss:2006}. While a fair number of traditional areas of mathematics are related or dedicated to lattice structures and periodicity, from crystallographic groups and quadratic forms to the geometry of numbers and sphere packings \cite{hilbert:cohn-vossen:GeometryImagination:1962,weyl:symmetry:1952,conway:sloane:SpherePackings:1999}, mathematical studies of periodic frameworks are of relatively recent date \cite{kotani:sunada:crystalLattice:2001,delgadoFriedrichs:PeriodicGraphs:2005,borcea:streinu:periodicFlexibility:2010}. 

\noindent
{\bf Displacive phase transitions.} Crystallographers have long observed that 
under variations of temperature or pressure, materials may undergo structural changes while retaining their crystalline nature. Understanding these phase transitions remains a  challenging problem \cite{dolino:quartzAB:1990}. The role of geometry is particularly relevant in {\em displacive} phase transitions which may be interpreted as continuous deformations of a given framework \cite{dove:displacive:1997}. A rigorous mathematical understanding of rigidity and flexibility
properties of periodic frameworks may assist, clarify and complement approaches based on computational physics and materials science
\cite{kapko:treacy:thorpe:guest:collapse:2009,borcea:streinu:periodicFlexibility:2010}.

\medskip
\noindent
{\bf Computational challenges.}
Flexibility studies of macromolecules modeled as mechanical (finite or periodic) frameworks \cite{goodwin:rums:2006,kapko:treacy:thorpe:guest:collapse:2009} are oftentimes conducted with computationally expensive yet numerically imprecise simulations. For finite structures characterized by theorems of Maxwell-Laman type (described below), much faster approaches based on degree-of-freedom counting and rigid component calculations are known. But such theorems are exceedingly rare and difficult to obtain. For periodic structures,  one would expect, intuitively, that the number and distribution of connections in a ``unit cell'' would dictate the flexibility properties of the framework, but the occurrence of certain `paradoxes' \cite{kapko:treacy:thorpe:guest:collapse:2009} related to choices of periodicity has stymied the progress. Making precise such an intuition with an appropriate mathematical definition has emerged only very recently \cite{borcea:streinu:periodicFlexibility:2010}, opening the way to a combinatorial treatment and efficient algorithms \cite{borcea:streinu:minimalPeriodic:LMS:2011} for rigidity and flexibility analysis of periodic {\em bar-and-joint} structures (made from rigid bars connected through rotatable joints).

\medskip

The  approach adopted in this paper is rooted in the theory developed in these recent papers \cite{borcea:streinu:periodicFlexibility:2010,borcea:streinu:minimalPeriodic:LMS:2011}, where 
we gave a Maxwell-Laman theorem for {\em generic} periodic bar-and-joint rigidity, thus reducing both the long-standing finite case and several other periodic situations to genericity refinement conjectures. It is therefore important to understand in this context the intrinsic theoretical difficulties we were facing in proving the results of this paper: the rest of the introduction gives this perspective.

\medskip
\noindent
{\bf Modeling crystalline frameworks.} Because of interatomic bond length and angle constraints, some crystalline materials 
may be abstracted as a periodic framework made from rigid bars connected at vertices, to which, in principle, the theory developed in \cite{borcea:streinu:periodicFlexibility:2010,borcea:streinu:minimalPeriodic:LMS:2011} applies. Often, substructures making individual {\em rigid bodies} are identified from the outset, such as the octahedra in Fig.~\ref{fig:prussianBlue} or, in a 2D example, the gray bodies from Fig.~\ref{fig:periodicFramework}(b). In this paper we focus on those made of rigid bodies connected by rigid bars, as illustrated in Fig.~\ref{fig:prussianBlue} and ~\ref{fig:periodicFramework}(b).

\begin{figure}[h]
 \centering
 \subfigure[]{\includegraphics[width=0.28\textwidth]{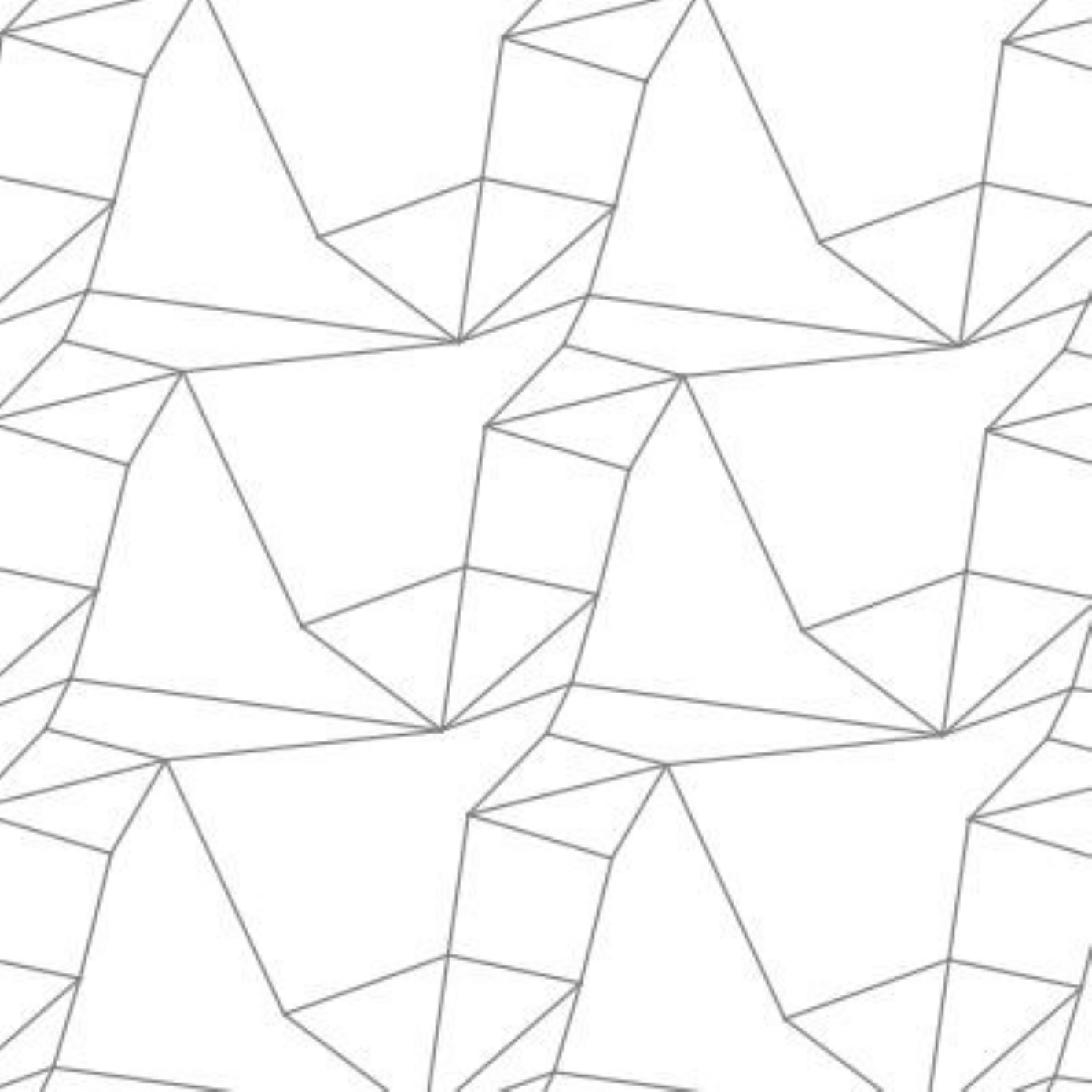}} \ \ \ 
 \subfigure[]{\includegraphics[width=0.28\textwidth]{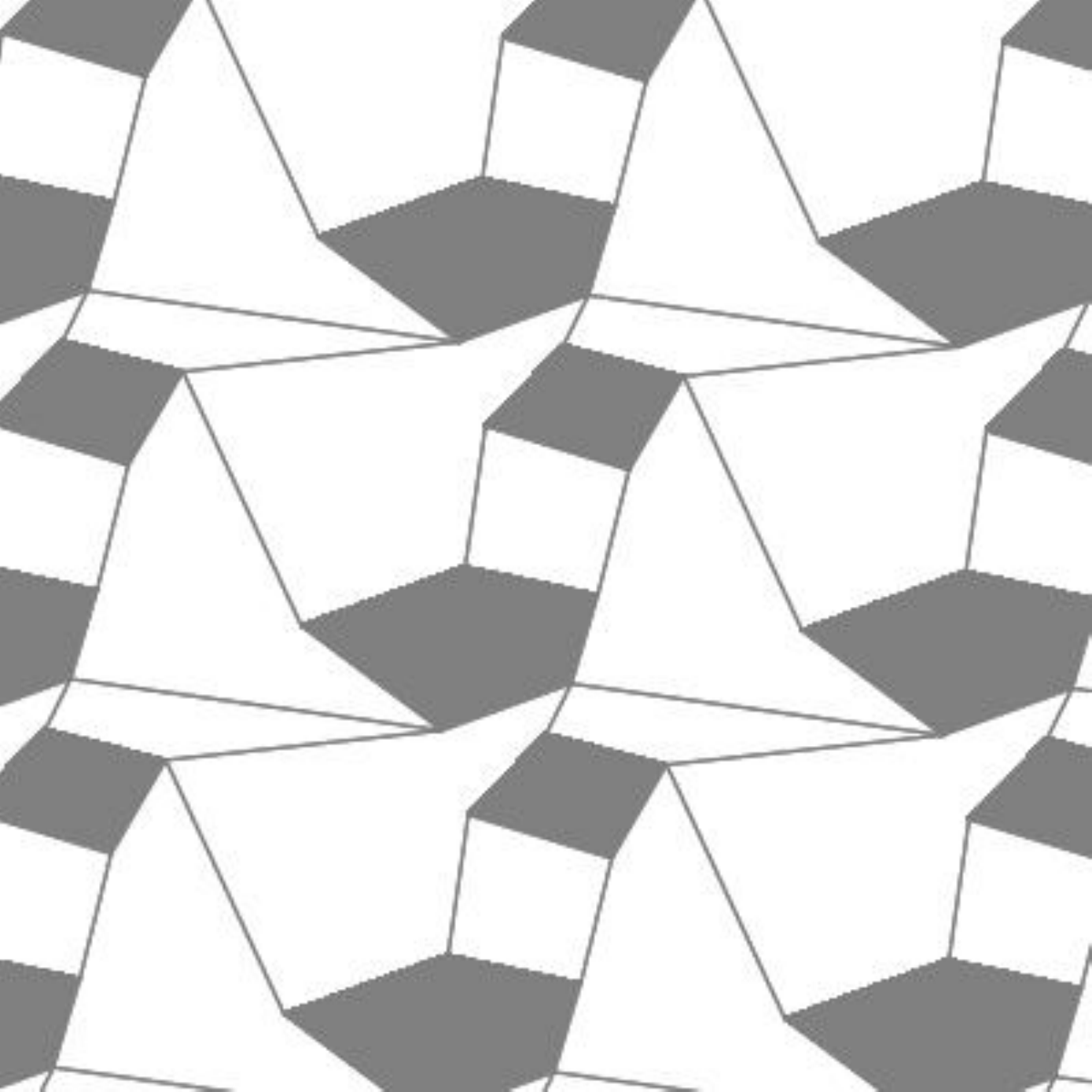}}
 \caption{ \small{(a) A 2D periodic bar-and-joint framework containing smaller rigid components. In (b), the same framework, viewed as a body-and-bar framework: the rigid components form the bodies, and the remaining bars connect distinct bodies.  }
}
 \label{fig:periodicFramework}
\end{figure}

The practical difficulty in applying the results of  \cite{borcea:streinu:minimalPeriodic:LMS:2011} to these frameworks comes from the {\em generic} nature of this characterization theorem. Indeed, the placement of the ``bars'' in actual crystal models, with several bars sharing endpoints and connecting pairs of atoms with specified regularity, is not guaranteed, {\em a priori}, to be generic.  This is a very subtle distinction for which the reader may need further guidance. Hence we discuss now briefly the challenges inherent in addressing questions of generic rigidity, at the same time presenting an overview of the techniques we use.

\medskip
\noindent
{\bf Maxwell-Laman sparsity conditions.} The theory of finite frameworks has a distinguished tradition going back almost 150 years to Maxwell \cite{maxwell:equilibrium:1864}, where a {\em sparsity condition} was shown to be necessary for minimal rigidity of bar-and-joint frameworks: in dimension $d$, for any subset of $d\leq n'\leq |V|$ vertices, the graph should span at most $dn'-{{d+1}\choose{2}}$ edges, with equality for the whole set of $n=|V|$ vertices. Its sufficiency for generic frameworks in dimension two was proven over 100 years later (Laman's theorem~\cite{laman:rigidity:1970}), and is known to fail in higher dimensions. The problem of completing the combinatorial characterization for bar-and-joint frameworks in arbitrary dimensions remains the most conspicuous open question in rigidity theory. However, some restricted classes of finite frameworks  have been shown to have similar Maxwell-Laman counting characterizations: body-and-bar and body-and-hinge frameworks~\cite{jackson:jordan:bodyBarHinge:EJC:2009,whiteley:unionMatroids:1988,tay:rigidityMultigraphs-II:1989,tay:rigidityMultigraphs-I:1984,tanigawa:dilworth:2010}, 
panel-and-hinge frameworks~\cite{katoh:tanigawa:proofMolecularConjecture:DCG:2011} and 
$(d-2)$-dimensional plate-and-bar frameworks~\cite{tay:linking:GraphsAndCombinatorics:1991,tay:rigidityMultigraphs-II:1989,tanigawa:dilworth:2010}.

\medskip
\noindent
{\bf Generic structures.} All known theorems giving combinatorial characterizations of framework rigidity and flexibility proceed by associating a {\em rigidity matrix $M$} to the framework. A set of algebraic constraints (e.g. fixed distances between some pairs of points) are applied on a collection of rigid objects (e.g. a finite set of points), with the goal of obtaining a rigid structure, i.e. one in which all the pairwise distances are determined. The tangent space at a generic point on the variety given by these algebraic constraints induces a matrix, called the {\em rigidity matrix}, whose rank we seek to compute. The maximal number of independent constraints that can be imposed is also the minimum number that would make the structure rigid. Together with additional constraints that eliminate the trivial rigid motions, they lead to a square rigidity matrix. 
Maximal rank, attained when the determinant $det(M)$ of the rigidity matrix is non-zero, indicates that the framework is {\em minimally rigid}: it is rigid and the removal of any constraint will make it flexible. The condition of maximal rank is therefore an algebraic condition in the free variables of the system. A framework (with a square rigidity matrix) is {\em generic} if its actual geometric realization falls outside the set of zeros of the algebraic equation $det(M)=0$. This extends naturally to frameworks whose rigidity matrix is not square, by expressing the rank condition in terms of corresponding minors of $M$. 

Theorems of Maxwell-Laman type for minimally rigid structures characterize classes of graphs underlying frameworks whose rigidity matrices have maximum rank, generically, for certain structural models as described above. A specific geometric realization of a graph belonging to such a class is, with probability $1$, rigid. But the measure-zero set when it is not rigid is non-empty, and {\em may} contain important cases that appear in practice. Refining the genericity conditions to prove that they include such practical occurrences is in general a very difficult problem. A remarkable and important success of this nature is the recent proof of the {\em Molecular Conjecture} \cite{katoh:tanigawa:proofMolecularConjecture:DCG:2011}.
Our contribution here is of a similar flavor, in the sense that it eliminates a certain type of {\em genericity assumption}. 

\medskip
\noindent
{\bf Maxwell-Laman sparsity for periodic frameworks.}
The connection pattern of a periodic body-and-bar system determines an infinite multi-graph $G=(V,E)$, with vertices corresponding to bodies and edges corresponding to bars between pairs of connected bodies. See Fig.~\ref{fig:periodicFramework}(c) for an example in 2D and Fig.~\ref{fig:prussianBlue} for one in 3D. Periodicity, or more precisely $d$-periodicity (where $d$ represents the dimension of the ambient space in which the graph is realized geometrically) requires a free Abelian automorphism subgroup $\Gamma\subset Aut(G)$ of rank $d$. We work under the assumption that the quotient graph $G/\Gamma$ has a finite number $n$ of vertex orbits and a finite number $m$ of edge orbits.

The problem of characterizing periodic frameworks is dramatically different from the finite case: the generic periodic bar-and-joint frameworks have been characterized in all dimensions \cite{borcea:streinu:minimalPeriodic:LMS:2011} by a Maxwell-sparsity condition {\em on the quotient graph}. The finite bar-and-joint frameworks, as well as the periodic body-and-bar, body-and-hinge, body-and-pin etc. appear as special cases where {\em additional algebraic dependencies} are present. As indicated, such cases are not guaranteed to be generic a priori, and
even getting a necessary sparsity condition  (something that was trivial in the finite case) is not easy.

\medskip
\noindent
{\bf Matroidal sparsity conditions.} An additional difficulty, relevant from an algorithmic point of view, arises from the fact that not all types of sparsity conditions are {\em matroidal}. For instance, Maxwell sparsity for finite 3D bar-and-joint frameworks (the ``3n-6'' condition) is not. For matroidal conditions, efficient recognition algorithms exist. Relevant for our setting are matroid partition \cite{edmonds:matroidPartition:1965} and pebble game algorithms \cite{streinu:lee:pebbleGames:2008}.  

We can now state our contributions in more technical terms.
  
\noindent
{\bf Contribution and novelty.} Motivated by foundational questions in theoretical crystallography, we present a complete mathematical and algorithmic solution to the {\em periodic body-and-bar problem} by giving a  Maxwell-Laman characterization. We formally define the periodicity model, identify the necessary sparsity conditions, prove their sufficiency and, finally, give efficient algorithms for deciding rigidity, based on matroid partitioning and pebble games.

\medskip
\noindent
As already seen in the bar-and-joint case \cite{borcea:streinu:periodicFlexibility:2010,borcea:streinu:minimalPeriodic:LMS:2011}, the mutual relations between the finite and
the periodic settings are intricate. We emphasize here that our approach 
considers the periodicity group $\Gamma$ as part of the initial data $(G,\Gamma)$ but
allows the variation of its representation as a  lattice of translations of full rank. Moreover,
the problem of minimal rigidity may be formulated at several levels. The most basic level
looks only at the structure of the quotient graph $G/\Gamma$ and provides characterizations `up to a generic lifting of edges' from $G/\Gamma$ to a covering periodic graph.
The main results presented below are formulated at this level. For $n=1$, that is, in the case
of a single body orbit under $\Gamma$, an answer is given at all levels and offers a first
measure of the contrast between solving in terms of $G/\Gamma$ or fully in terms of
the given periodic graph $(G,\Gamma)$. We discuss aspects implicated at this higher level 
after solving the problem at the basic level.

\medskip
\noindent
{\bf Related work.} Among other recent contributions to the subject, we mention \cite{malestein:theran:combRigidityPeriodic2D:arxiv:2010} which addresses the two-dimensional case of periodic bar-and-joint frameworks and \cite{ross:schulze:whiteley:finiteMotions:2011,power:crystal:arxiv:2011} where certain
classes with `added symmetry' are investigated. The final part of the thesis \cite{ross:thesis:2011}  includes
some considerations about body-and-bar frameworks ``on the fixed torus". This setup
corresponds with a fixed lattice of periods. A necessary condition for minimal rigidity is
presented and conjectured to be also sufficient.  It will be seen below that this is a restricted
version of the minimal rigidity problem in terms of $(G,\Gamma)$.

\medskip
\noindent
{\bf Overview.}
We describe our setting in Section~\ref{sec:gf}, where we also derive the rigidity matrix. We
indicate different levels for the problem of minimal rigidity in Section~\ref{sec:mr} and solve it at all levels for the base case $n=1$. Our main result characterizing minimal rigidity in terms of the structure of the quotient graph $G/\Gamma$ is proven
in Section~\ref{sec:mainMR}. We then extend
these results in Section~\ref{sec:mixed} to mixed periodic plate-and-bar frameworks. 
We conclude with the algorithmic consequences and several open problems.


\section{Periodic graphs and frameworks}
\label{sec:gf}

In this section we introduce our formal definition of periodic rigidity for periodic body-and-bar frameworks. For standard rigidity theoretic definitions such as infinitesimal rigidity and rigidity matrix, see \cite{borcea:streinu:periodicFlexibility:2010,whiteley:unionMatroids:1988}.

\begin{figure}[h]
 \centering
 \subfigure[]{\includegraphics[width=0.60\textwidth]{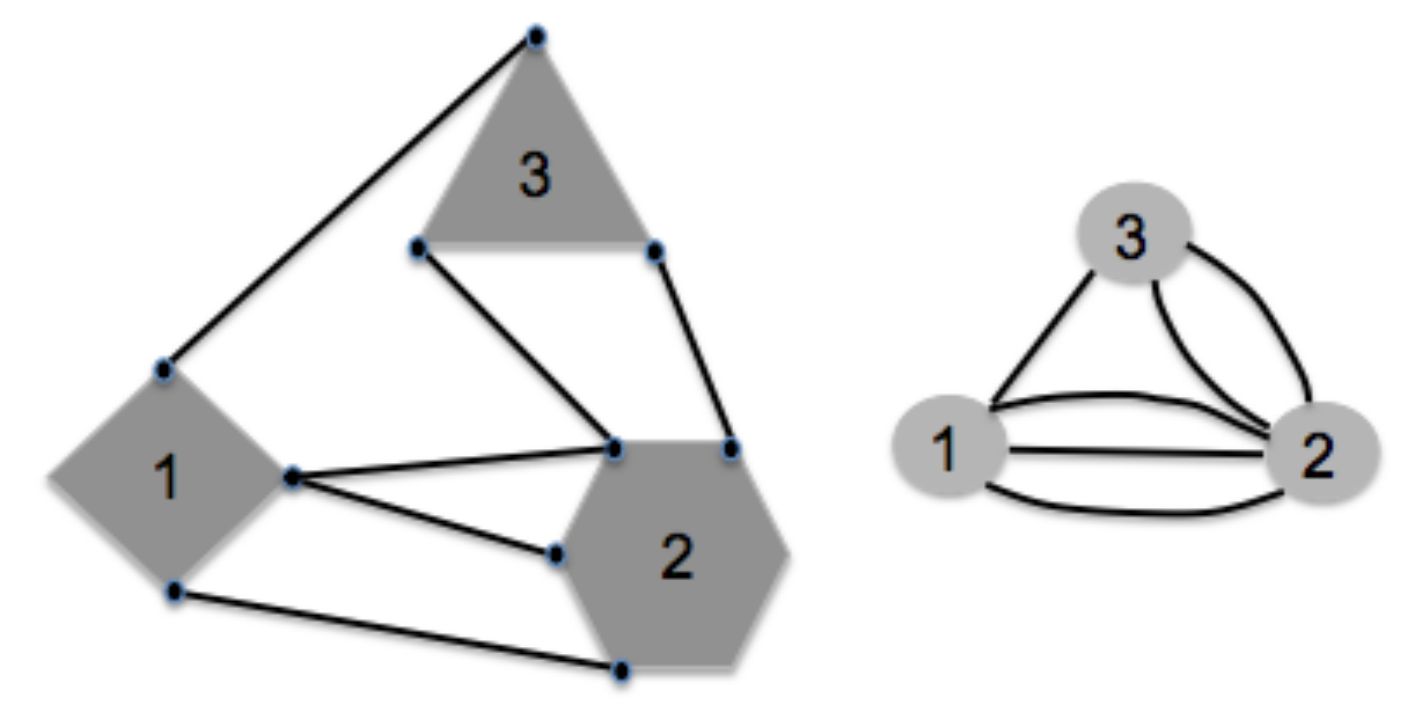}}
 \caption{\small{A finite, minimally rigid body-and-bar framework in 2D, with 3 bodies and 6 bars, and its associated multi-graph.}} 
 \label{fig:bodyBar2D}
\end{figure}

We remind the reader a classical concept from finite rigidity theory: a finite {\bf body-and-bar} $d$-framework is a finite collection of rigid bodies in $R^d$ connected through bars. The endpoints of a bar lie on two distinct bodies and act as rotatable joints, i.e. they permit rigid motions of the incident objects (bodies or bars) relative to each other, while contact is maintained at the joints. Several bars may be placed between the same pair of bodies. See Fig.~\ref{fig:bodyBar2D}. The combinatorial (incidence) structure of such a framework is captured by a multi-graph $G=(V,E)$, whose vertices $V$ correspond to bodies and edges to bars. 
Body-and-bar frameworks have been studied in several foundational papers (e.g. \cite{tay:rigidityMultigraphs-I:1984,tay:rigidityMultigraphs-II:1989}) and have practical application in studies of molecular flexibility \cite{thorpe:etAl:flexibilityBiomolecules:2005,fox:etAl:KINARIweb:2011}. 
We proceed now to our new definitions.

\medskip
\noindent
A {\bf $d$-periodic body-and-bar graph} is a pair $(G,\Gamma)$ of an infinite multi-graph $G=(V,E)$ and a group $\Gamma$ acting on it. 
The group $\Gamma$, called the {\em periodicity group}  
of $G$, is a rank $d$ free Abelian subgroup of the automorphism group $Aut(G)$ of $G$, which acts without fixed points and has a finite number of vertex and edge orbits. The elements $\gamma \in \Gamma$ may be called {\em periods} of $G$. 
Since both $n = card(V/\Gamma)$ and $m = card(E/\Gamma)$ are finite, this setting induces a finite quotient multi-graph $G/\Gamma=(V/\Gamma, E/\Gamma)$.

\noindent
{\bf Isometries.} We denote  by $E(d)$ the Euclidean group in dimension $d$, that is, the isometry group of the Euclidean space $R^d$. The subgroup of translations is denoted by ${\cal T}(R^d)\subset SE(d)$. The connected component of the identity, denoted by $SE(d)$, is made of all the orientation preserving isometries and is referred to as the special Euclidean
group or the group of rigid motions in $R^d$. A transformation in $E(d)$, respectively $SE(d)$,
can be represented by a pair $T=(p,M)$, where $ p\in R^d$ denotes a translation and $M\in O(d),\ \mbox{respectively}\ M\in SO(d)$, is an orthogonal transformation, resp. special orthogonal one. $M$ is a matrix with $|det( M )|=1$, resp. $det( M ) = 1$. 
An Euclidean transformation $T: R^d \rightarrow R^d$ operates by the formula $T(x)=Mx+p$. In this notation translations correspond to pairs $(p,I_d)$, where $I_d$ is the $d\times d$ identity matrix.

\begin{figure}[h]
 \centering
 \subfigure[]{\includegraphics[width=0.60\textwidth]{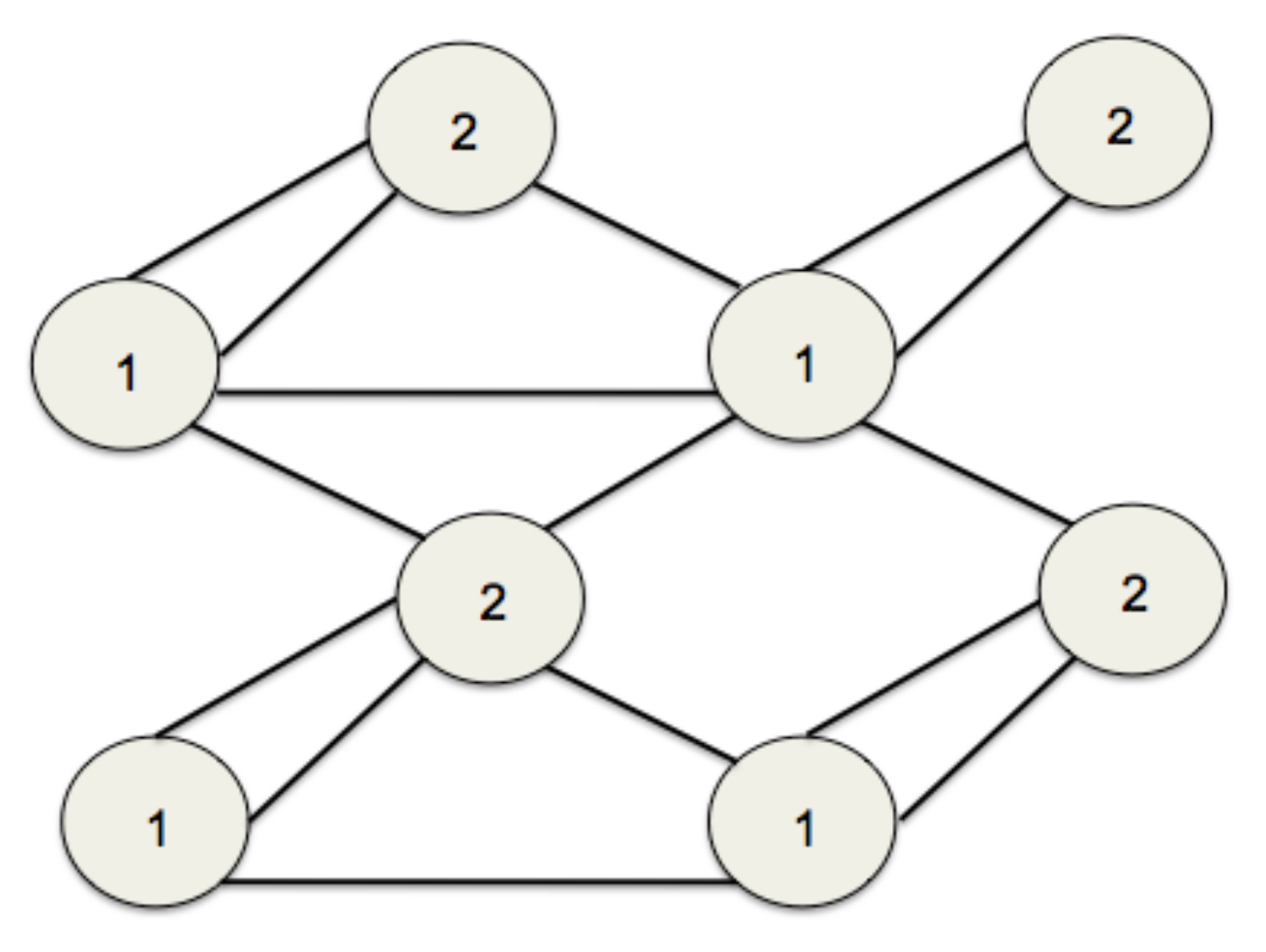}}
    \caption{\small{Top: A small fragment of an abstract periodic graph. Only the edges between depicted vertices are drawn. Equivalent vertices under the action of the group of periods have the same label. The quotient multi-graph and a realization are shown in Fig.~\ref{fig:quotientGraph} and Fig.~\ref{fig:framework}. 
}} 
    \label{fig:periodicGraph}
\end{figure}

\medskip
\noindent
To define a body-and-bar framework, we will specify a placement in $R^d$ of the bodies. We assume that all our bodies are oriented as the ambient space $R^d$. Then the placement of a body in $R^d$ is described by a pair $T=(p,M)\in SE(d)$. In other words, we conceive of the body as described by a (positively oriented) Cartesian frame with origin at $p$ and basis
vectors corresponding to the columns of the orthogonal matrix $M\in SO(d)$.

Our next goal is to define a {\bf periodic body-and-bar framework}. A 2D example  is given in Fig.~\ref{fig:framework}. 

\medskip
\noindent
A {\bf presentation} in $R^d$ of a $d$-periodic body-and-bar graph $(G,\Gamma)$ 
is given by an assignment of Cartesian frames to the vertices $\tau: V \rightarrow SE(d)$, alongside with an injective 
representation of the periodicity group $\Gamma$, $\pi: \Gamma \rightarrow {\cal T}(R^d)$. We refer to $\pi(\Gamma)$ as the {\em lattice of periods} and ask that it has rank $d$.

\medskip 
\noindent
Furthermore, for each edge $e=(i,j)\in E$ two (arbitrary but fixed) endpoints are given, indicating where the bar is attached on the bodies corresponding to vertices $i$ and $j$. Let us denote these endpoints with $q^i=q^i(e)$, $q^j=q^j(e)$. 
The coordinates $q^i\in R^d$ are given with respect to the Cartesian frame marking the body corresponding to $i$, so that for a frame
$\tau(i)=(p_i,M_i)$, this end of the bar is at the point of $R^d$ with coordinates $  \tilde{q}^i=M_iq^i+p_i $. 

\medskip 
\noindent
All this data must respect the conditions required by periodicity, i.e. when we act on a vertex $i$ by a period $\gamma\in \Gamma$, we have $\tau(\gamma i)=\pi(\gamma) \tau(i)$ and $q^i(e)=q^{\gamma i}(\gamma e)\ \ \mbox{etc.}$ The main definition can now be given:

\medskip 

\begin{figure}[h]
 \centering
 \subfigure[]{\includegraphics[width=0.30\textwidth]{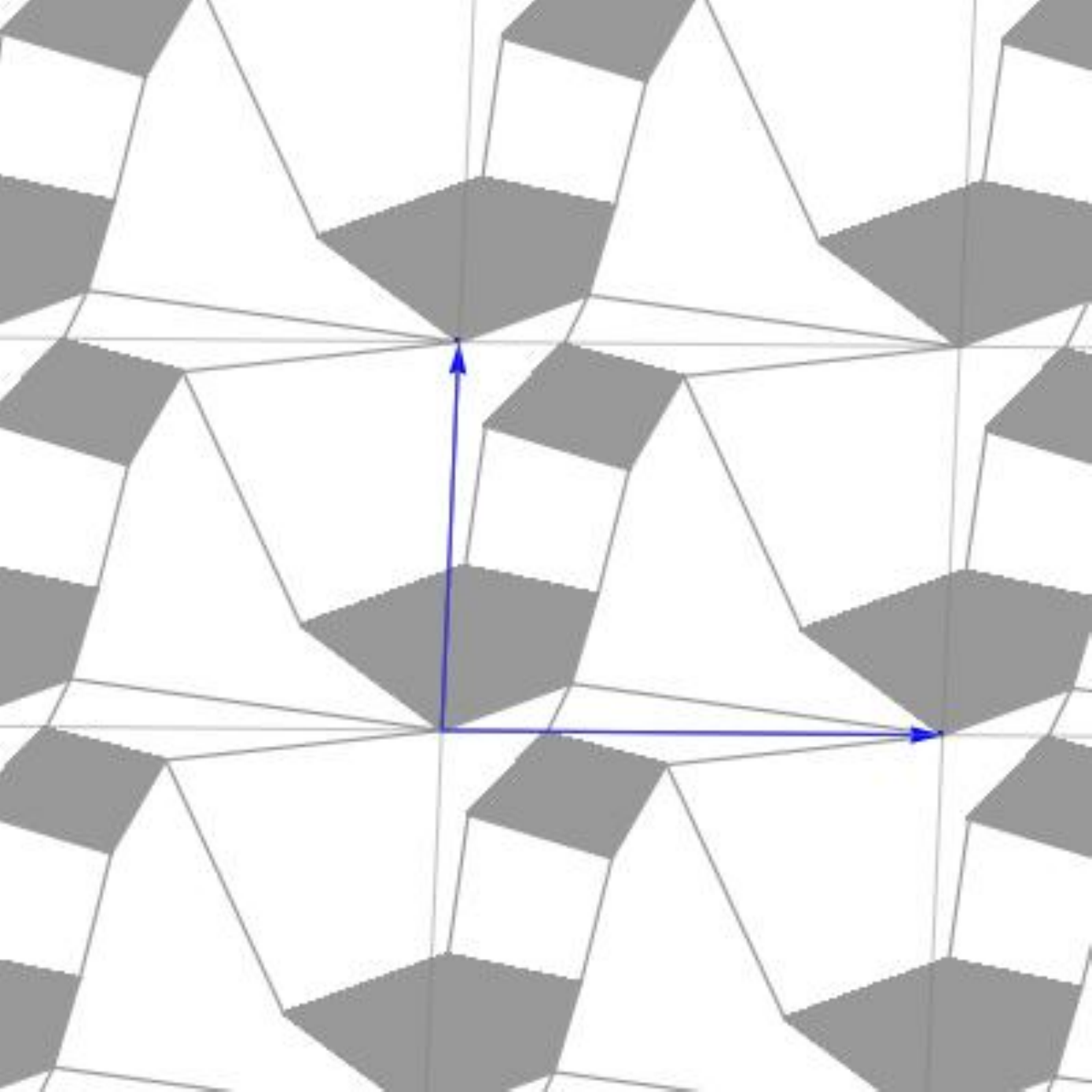}} \ \ \ 
 \subfigure[]{\includegraphics[width=0.30\textwidth]{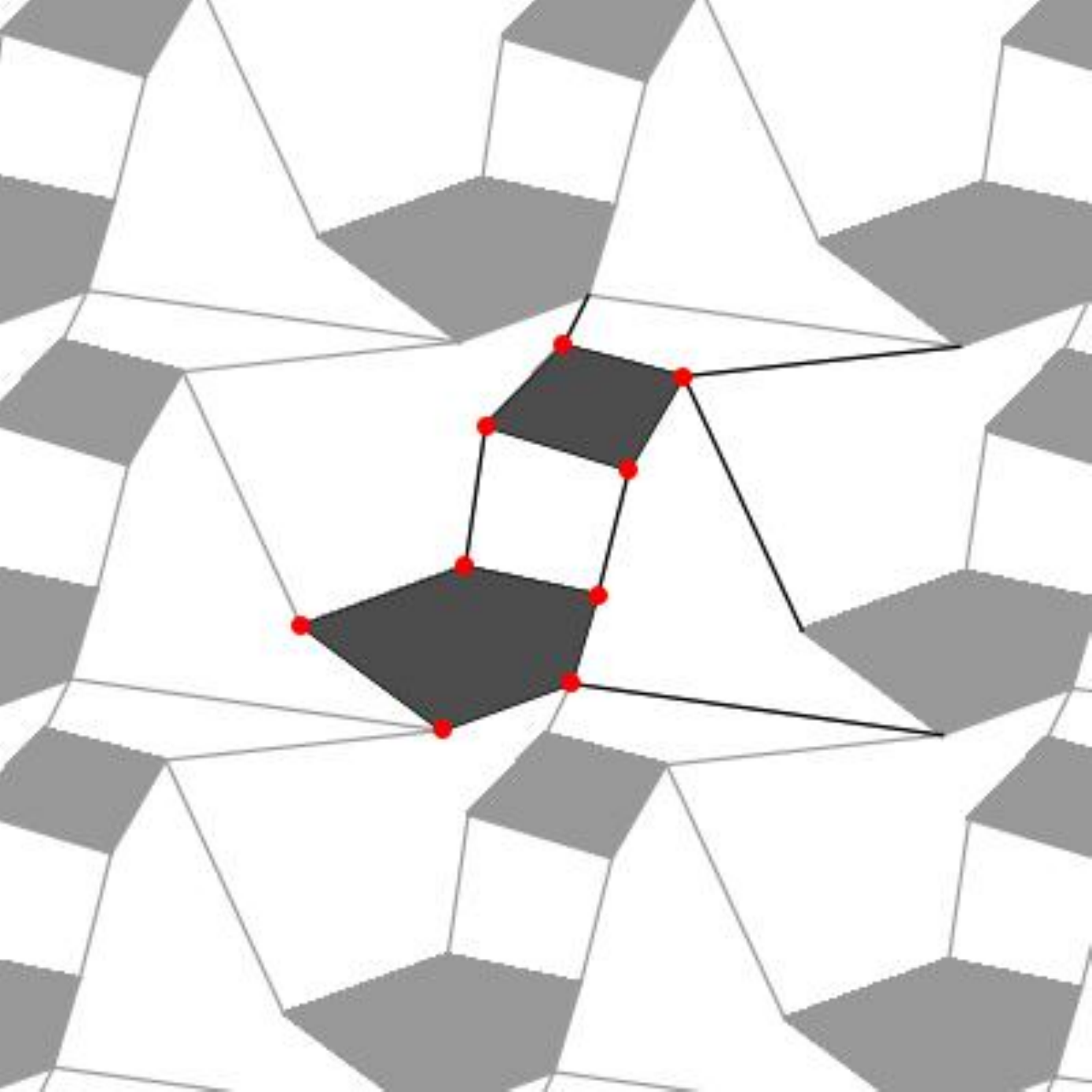}} \ \ \ 
 \subfigure[]{\includegraphics[width=0.30\textwidth]{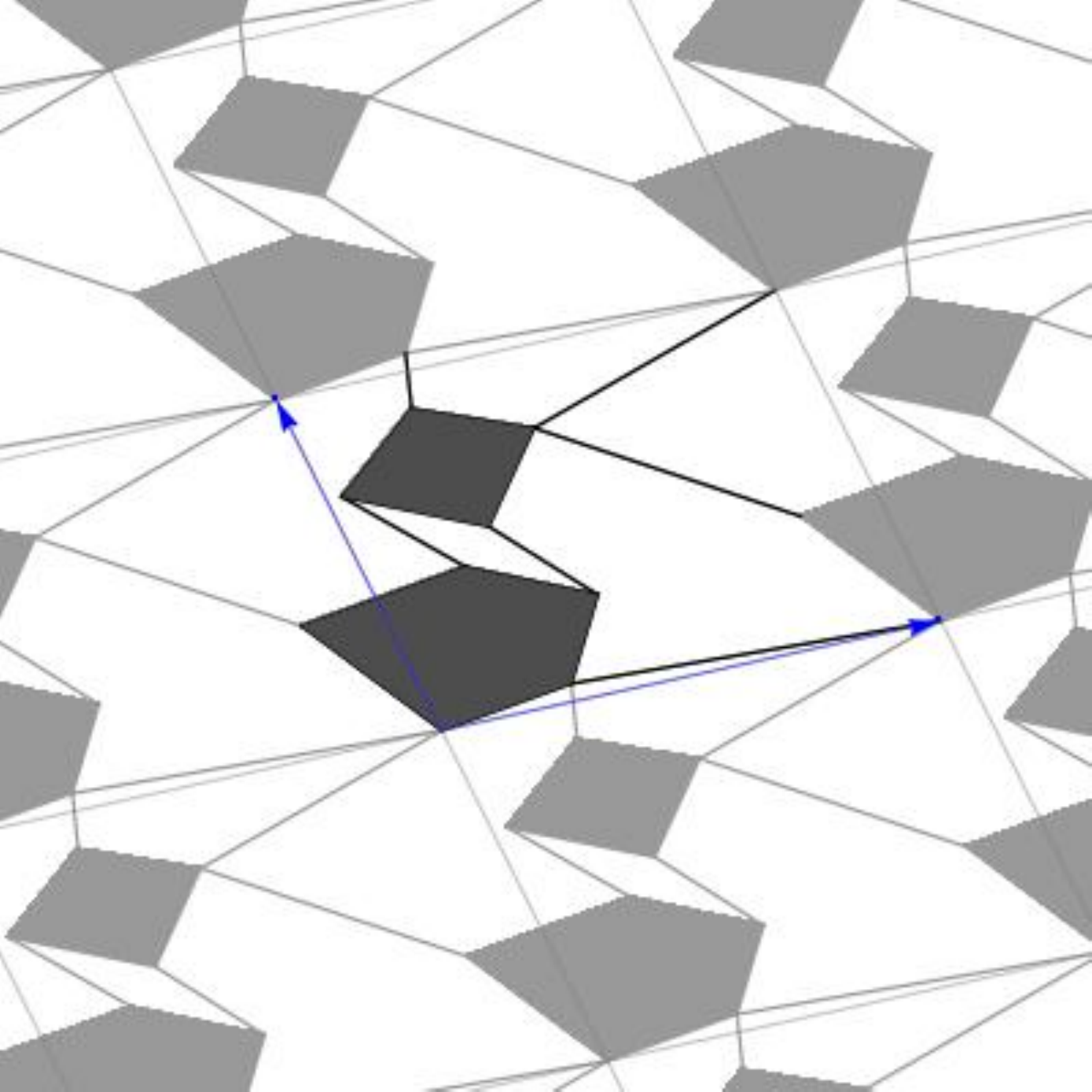}}
\vspace{-0.15in}
 \caption{\small{ (a) Fragment of a $2$-periodic body-and-bar framework ($d=2$), with the generators (which are not bars) of the periodicity lattice marked by arrows and the orbit of one bar endpoint visualized as a ``lattice''.
(b) Representatives of the $n=2$ equivalence classes of bodies and $m=6$ classes of bars are shown in darker shades. The endpoints of the bars, marked in red, are not all distinct in this case.
(c) The framework is {\em periodically flexible}, as illustrated by this deformation of (b). Notice the deformation of both the framework and of the periodicity lattice.
}}
 \label{fig:framework}
\vspace{-0.2in}
\end{figure}

\medskip
\noindent
A $d$-{\bf periodic body-and-bar framework} is a connected $d$-periodic graph $(G,\Gamma)$
together with a presentation $(\tau,\pi,q)$ where all bars have non-zero lengths.

\medskip
\noindent {\bf Realizations, configurations, rigid and flexible frameworks.}
Presentations which use the same attachment points $q$ (relative to the corresponding bodies) and induce the
same system $\ell$ of bar lengths form the {\em realization space} for the data $(G,\Gamma, q, \ell)$.

Realizations that differ by an isometry of $R^d$ will be considered as the same {\em configuration}, and the
{\em configuration space} of $(G,\Gamma,\tau,\pi,q)$ is defined as the {\em quotient space} of all realizations by the group
$E(d)$ of all isometries of $R^d$.  The deformation space of a periodic framework $(G,\Gamma,\tau,\pi,q)$ is the
connected component of the corresponding configuration. When the deformation space consists of a single point,
the framework is called {\em rigid}, otherwise it is {\em flexible}.
An example of a flexible periodic framework is illustrated in Fig.~\ref{fig:framework}, with its quotient graph shown in Fig.~\ref{fig:quotientGraph}. %

\begin{figure}[h]
 \centering
 \subfigure[]{\includegraphics[width=0.40\textwidth]{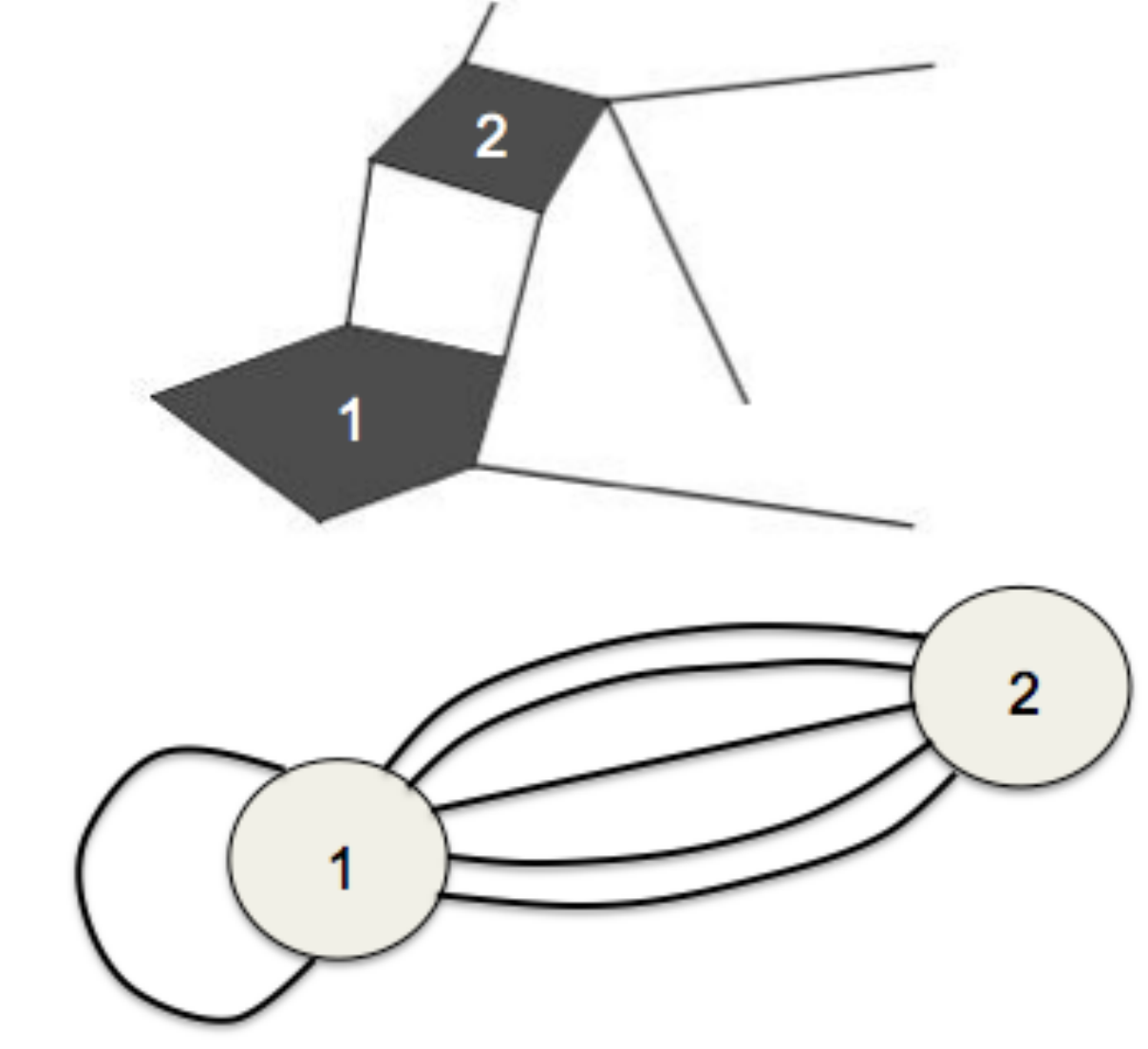}}
    \caption{\small{Top: representatives of orbits for bodies and bars for the periodic framework from Fig.~\ref{fig:framework}. Bottom: the quotient multi-graph. 
}} 
    \label{fig:quotientGraph}
\end{figure}

\medskip
\noindent
As in \cite{borcea:streinu:periodicFlexibility:2010}, this deformation theory has the following useful characteristics:
{\bf(i)} The periodicity group $\Gamma$ is part of the structure and the deformations of a framework are those
preserving this specified periodicity.
{\bf(ii)} The lattice of translations $\pi(\Gamma)$ representing $\Gamma$ is allowed to vary as the framework deforms.
{\bf(iii)} The realization space is the solution space of a finite algebraic system of quadratic equations.
{\bf(iv)} There is an infinitesimal deformation theory, obtained by differentiating the equations. This leads to the concept of {\em infinitesimally rigid framework}, characterized by a {\em rigidity matrix} of maximum rank.

\medskip
\noindent
In particular, standard arguments from algebraic geometry (inverse function theorem) can now be used to prove that infinitesimal rigidity implies rigidity.

\medskip
\noindent
{\bf Minimal rigidity.} Following a heuristic going back to Maxwell \cite{maxwell:equilibrium:1864}, we do a quick calculation of the minimal number of parameters needed to specify a rigid framework. This gives us an indication of the number of bars we would expect in a minimally rigid graph. A vertex orbit requires ${{d+1}\choose{2}}$ parameters to be specified, and an additional $d^2$ parameters specify the periodicity lattice. The total number of variables in the algebraic system is thus $n{{d+1}\choose{2}}+d^2$. Each edge orbit leads to a constraint (equation), and generically it eliminates one degree of freedom. Since we cannot eliminate the trivial isometries, we expect to have at most $n{{d+1}\choose{2}}+d^2-{{d+1}\choose{2}}=(n-1){{d+1}\choose{2}}+d^2$ independent equations. By definition:

A {\bf minimally rigid $d$-periodic body-and-bar graph} $(G,\Gamma)$ is a $d$-periodic graph with $ m = (n-1){{d+1}\choose{2}}+d^2$ which has an infinitesimally rigid presentation.

For $n=2$ bodies in dimension $2$ this number is $7$. Since the quotient graph from Fig.~\ref{fig:quotientGraph} of the periodic framework from Fig.~\ref{fig:framework} has $6<7$ bars, we expect the periodic framework to be flexible with one degree of freedom (1dof).

We have now all the ingredients to start developing the rigidity theory of periodic body-and-bar frameworks. The first step is to define and analyze the rigidity matrix associated to a generic framework.


\medskip
\noindent
{\bf The rigidity matrix.}
We express the rigidity matrix in terms of the following coordinates:  we choose a basis for the lattice of periods $\pi(\Gamma)$: $\lambda_1,...,\lambda_d$ and denote by $\Lambda$ the
$d\times d$ matrix with columns $\lambda_i$; then we choose representatives for the
$n$ equivalence classes of vertices modulo $\Gamma$ and consider the corresponding bodies  as marked by
$(p_i,M_i),\ i=1,...,n$. We then choose an orientation of the edge set in the quotient graph $G/\Gamma$ and select representatives for the equivalence classes of edges modulo $\Gamma$ by the rule 
that all edge representatives originate at the already chosen vertex representatives.  
We remind the reader that the endpoints $q^u$, $q^v$ of an edge representative are fixed parameters, but their choice is only subject to $\tilde{q}^u\neq\tilde{q}^v$. 
In order to keep the notation as simple as possible, we do not pursue the details related to marking distinctions
between several bars possibly joining the same pair of bodies.

\medskip
\noindent
Since one end (the origin) of any edge representative is on a body representative, the other end is on the corresponding representative translated by a period $\lambda=\Lambda C(e) =\Lambda c$, where $c=C(e)$ is a column vector with integer entries $c_i, \ i=1,...,d$.

Then, a bar constraint corresponding to an edge $e$ between bodies $i$ and $j$, with endpoints $q^i(e)$, resp. $q^j(e)$, on these two bodies, takes the form $|| \Lambda C(e) +M_jq^j(e)-M_iq^i(e) +p_j-p_i ||^2= \ell_e^2$, with $e$ running over all $m=|E/\Gamma|$ edge representatives.

\medskip
\noindent
We will eliminate the equivalence under rigid motions by assuming that the first representative
body is fixed, with $(p_1,M_1)=(0,I_d)$.
We note here that for an orthogonal matrix $M$, an `infinitesimal variation' $\dot{M}$ takes the
form $\dot{M}=MA$, with  $A$ skew-symmetric. Thus, if we differentiate the system   at $(p_2,M_2),...,(p_{n},M_{n}),\Lambda$, we obtain a linear
system in infinitesimal variations $(\dot{p}_i,A_i),\ i=2,...,n$ and $\dot{\Lambda}$ with rows
\begin{equation}\label{eq:infinitesimalConstraint}
\langle \dot{\Lambda}c +M_jA_jq^j-M_iA_iq^i+\dot{p}_j-\dot{p}_i, \Lambda c +M_jq^j-M_iq^i+p_j-p_i \rangle =0
\end{equation}
\noindent
where we have dropped the notational dependency on $e$ and used indices $i,j$ for the chosen
body-representatives, with $M_1=I_d$, $p_1=0$ and $A_1=0$. The {\bf rigidity matrix} is the matrix of this linear system (for some specific ordering of the unknowns).

\medskip 
\noindent
{\bf Normalized coordinates.}\ Given any $d$-periodic body-and-bar framework, the manner
of marking chosen representatives in each vertex orbit by Cartesian frames $(p_i,M_i)$ is also
a matter of choice and we may adopt the {\em normalization} which has all $p_i=0$ and $M_i=I_d$. The rigidity matrix is then the matrix of the linear system in  $\dot{p}_i$, 
$A_i=-A_i^t$ and $\dot{\Lambda}$:
\begin{equation}\label{eq:normalized}
 \langle \dot{\Lambda}c +A_jq^j-A_iq^i+\dot{p}_j-\dot{p}_i, \Lambda c +q^j-q^i\rangle =0
\end{equation}
\noindent
where the edge vectors $\Lambda c +q^j-q^i$ run over a complete set of edge representatives.


\section{Minimal rigidity}
\label{sec:mr}

As already observed, up to Euclidean rigid motions, a periodic body-and-bar framework is
described by ${d+1\choose 2}(n-1)+d^2$ parameters. A bar equivalence class
imposes a single constraint on these parameters hence rigidity requires at least  
${d+1\choose 2}(n-1)+d^2$ constraints.  A $d$-periodic body-and-bar graph $(G,\Gamma)$, with $n=|E/\Gamma|$ and $m=|E/\Gamma|={d+1\choose 2}(n-1)+d^2$, is called minimally rigid if it can be presented as an infinitesimally rigid framework in $R^d$, that is, as a framework with rigidity matrix of maximal rank $m$. 

\medskip
\noindent
{\bf Liftings and marked quotient graphs.} The problem of characterizing the structure of minimally rigid periodic graphs can be 
formulated at {\em three distinct levels}, illustrated in Fig.~\ref{fig:labeledQuotient}. The most discerning and demanding level concerns the periodic graph $(G,\Gamma)$. Two other less discerning  levels focus only on the quotient multi-graph $G/\Gamma$, with or without some additional {\em lifting} information retained.  

\medskip
\noindent
The information lost upon passing from the periodic framework $(G,\Gamma)$ to the quotient graph $G/\Gamma$ is the {\em lifting} of the edges to specific vertex representatives of the same orbit. This information can be retained in the following form. We fix a basis of $\Gamma$, that is, an isomorphism $\Gamma\approx Z^d$. Then we choose vertex representatives and an orientation of the quotient graph. By using edge representatives which originate at the chosen vertex representatives, we have a well defined {\em lifting function} $ C :  E/\Gamma \rightarrow Z^d $ which allows retrieval of $(G,\Gamma)$ from $G/\Gamma$. When several bars connect the same bodies, we can retain this information by joining the vertex representatives with single edge, marked with a number indicating the edge multiplicity, as in Fig.~\ref{fig:labeledQuotient}(b). This `multi-graph with multiplicities' is denoted by ${\cal G}/\Gamma$. Combining the edge multiplicities with the lifting function yields a labeled multi-graph with multiplicities, such as illustrated in Fig.~\ref{fig:labeledQuotient}(c). This contains all the information needed to reconstruct $(G,\Gamma)$.

\begin{figure}[h]
 \centering
 \subfigure[]{\includegraphics[width=0.28\textwidth]{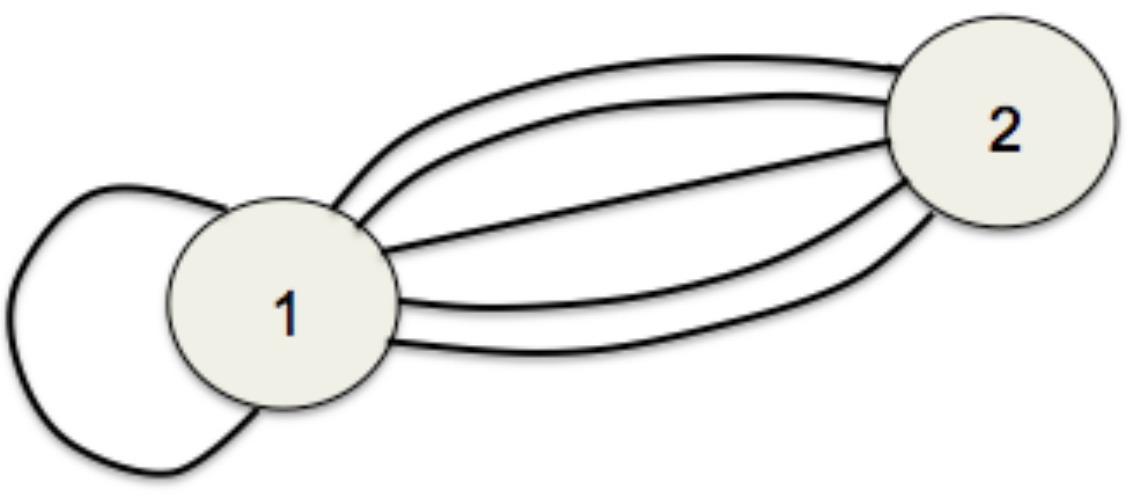}} \ \ \ 
 \subfigure[]{\includegraphics[width=0.28\textwidth]{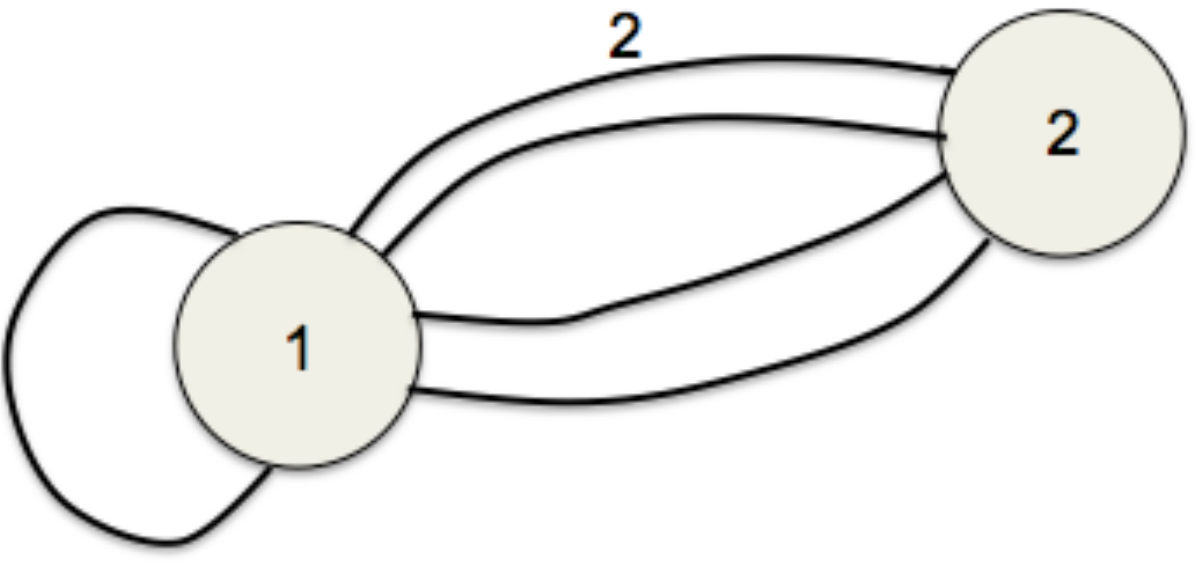}} \ \ \ 
 \subfigure[]{\includegraphics[width=0.34\textwidth]{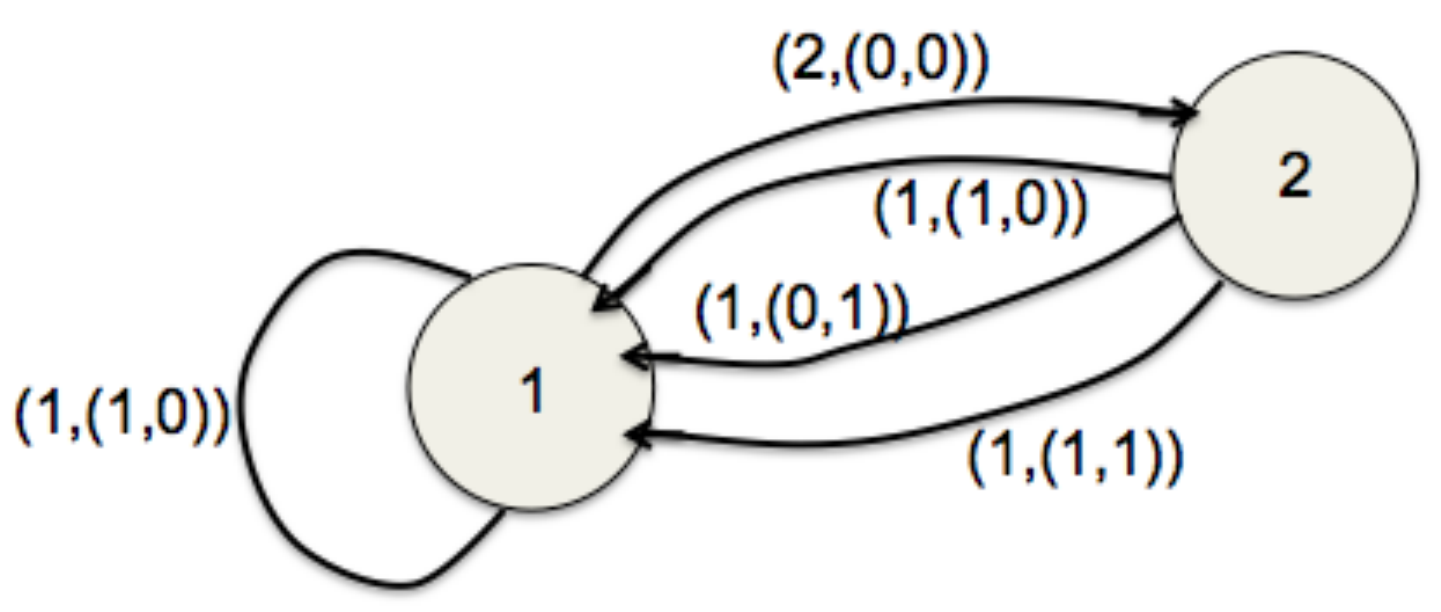}}
 \caption{ (a) The basic quotient graph $G/\Gamma$.   
(b) The quotient graph with multiplicities ${\cal G}/\Gamma$.
(c) The labeled and oriented quotient multi-graph with multiplicities, equivalent with $(G,\Gamma)$.
}
 \label{fig:labeledQuotient}
\end{figure}

\medskip 
\noindent
{\bf Remark:}\ An {\em invariant} way for describing the lifting data can be given in terms of 
the first {\em homology group} (with integer coefficients) $H_1(G/\Gamma)$, where $G/\Gamma$ is now seen as the topological space corresponding to the one-dimensional CW-complex defined by the quotient graph \cite{kotani:sunada:crystalLattice:2001}. With $G$ considered in the same topological perspective, we have a (connected) Abelian covering $G\rightarrow G/\Gamma$ with covering group $\Gamma$, and this is equivalent to a {\em surjective group homomorphism} $H_1(G/\Gamma) \rightarrow \Gamma$.

\medskip
\noindent
In the rest of this paper, we give a complete characterization of $(G,\Gamma)$ for the case of a single body $n=1$, and a  characterization in terms of the quotient graph $G/\Gamma$ for arbitrary $n$. 


\medskip 
\noindent
{\bf The $n=1$ case.}
The characterization of minimal rigidity in terms of $G/\Gamma$ is trivial for $n=1$, since the only
condition about the quotient graph is that of having $d^2$ loops. This is enough for finding generic liftings which have infinitesimally rigid framework presentations. We address here the more refined questions related to $(G,\Gamma)$ and ${\cal G}/\Gamma$, assuming $m=d^2$.

\medskip
The quotient multi-graph with multiplicities ${\cal G}/\Gamma$ has a single vertex and 
a number of loops with multiplicities $k_i,\ i=1,...,\ell$ with $ \sum_{i=1}^{\ell} k_i=d^2 $.

\medskip 
Minimal rigidity in this setting is  equivalent with the condition that all multiplicities
$k_i$ be at most $d$. This will follow from the solution of the sharper problem about
the structure of $(G,\Gamma)$.

\medskip
Recalling that $n=1,m=d^2$ and that we eliminate the equivalence under rigid motions
by fixing $(p_1,M_1)=(0,I_d)$, we obtain from (\ref{eq:infinitesimalConstraint}) a rigidity
matrix with $d^2$ rows of the form $C(e) \otimes [\Lambda C(e) +(q^j(e)-q^i(e))]$, where the tensor notation stand for a listing of all the $d^2$ products of the $d$ components
of the two indicated vectors. Since the endpoints of bars can be chosen at will, the edge vectors
$\Lambda C(e) +(q^v(e)-q^u(e))$ are arbitrary (non-zero) vectors, say $h=h(e)$. With simplified
notation $c_i$ for the coordinates in $C(e)$, such a row can be presented as $[\ h_1(c_1\ c_2\ ...\ c_d) \ h_2(c_1\ c_2\ ...\ c_d) \ ... \ h_d(c_1\ c_2\ ...\ c_d)\ ]$.

\medskip 
In this form one recognize a union operation on $d$ copies of a linear matroid,
leading to

\begin{theorem}
\label{thm:oneBodyOrbit}
Let $(G,\Gamma)$ be a $d$-periodic body-and-bar graph with $n=|V/\Gamma|=1$
and $m=|E/\Gamma|=d^2$. Then $(G,\Gamma)$ has infinitesimally rigid framework
presentations in $R^d$ if and only if for any set of equivalence classes of edges
$F\subset E/\Gamma$ we have 

\begin{equation}
\label{eq:minimal}
|F|\leq d\cdot dim( span C(F))
\end{equation}

\end{theorem}

\medskip 
\noindent
{\em Proof:}\ The necessity follows directly from the form, described above, 
of the row constraints. The sufficiency follows from results in matroid theory, which will be
used again in arguments for subsequent statements. Specifically, $d\cdot dim( span C(F))$ is an increasing submodular function
with values in the non-negative integers \cite{oxley:matroid:1992}. By \cite{pym:perfect:submodular:1970}, the associated matroid
is the union of $d$ copies of the matroid associated to the function  $dim(span C(F))$. Hence
a basis of the former is a disjoint union of bases in the latter which yield the maximal rank
of the rigidity matrix.\qed

\medskip 
\noindent
This proof also implies our earlier claim that a quotient graph ${\cal G}/\Gamma$ comes from
a minimally rigid $({\cal G},\Gamma)$ if and only if all multiplicities $k_i,\ i=1,...,\ell$ are at most $d$. 

Indeed, any choice of lifting function  $ C : E/\Gamma \rightarrow Z^d$, with
image consisting of classes of at most $d$ vectors (corresponding to the given multiplicities) and with representatives in general position, will satisfy the criterion of  Theorem~\ref{thm:oneBodyOrbit}.


\section{Minimal rigidity characterization of the quotient graph $G/\Gamma$}
\label{sec:mainMR}

In this section we present the main theoretical result. We assume that $m={d+1\choose 2}(n-1)+d^2$ and characterize the 
quotient graphs $G/\Gamma$ which are obtained from minimally rigid $d$-periodic
body-and-bar graphs $(G,\Gamma)$. In other words, we characterize minimal rigidity
`up to a generic lifting of the edge representatives'.

\medskip
\noindent
{\bf Graph sparsity.} We say that a multi-graph $G=(V,E)$ (possibly with loops) is $(a,b)$-sparse, or has {\em $a n-b$ sparsity type}, if any subset of $n'\leq |V|$ vertices spans at most $a n'-b$ edges. When equality holds we have a {\em tight} (sub)graph. See \cite{streinu:lee:pebbleGames:2008} for a comprehensive treatment of this kind of graph sparsity.

When $0\leq b<2a$, in particular when $b=0$ and $b=a$, the family of $(a,b)$-sparse graphs on a fixed number of vertices $n$ form a matroid. When $b=a$, this matroid is the union of $b$ graphic (tree) matroids. When $b=0$, it is the union of {\bf $a$} bi-cycle matroids, whose bases are the $(1,0)$-sparse graphs (also known in the literature as  map-graphs, pseudo-forests or bi-cycles). We also obtain a matroid when $b<0$: as shown in Corollary \ref{cor:anbSparse} below, it is the union of the $(a,0)$-sparsity matroid and the uniform matroid on $b$ elements, or, equivalently, it is any graph which has $an+b$ edges and contains a $(a,0)$-sparse subgraph. 

We characterize now the quotient graphs of periodic minimally rigid graphs in terms of matroid unions of graph sparsity matroids.

\begin{theorem}
\label{thm:minRigidity}
A  quotient graph $G/\Gamma$  with $m={d+1\choose 2}(n-1)+d^2$ is the quotient graph of an infinitesimally rigid periodic body-and-bar framework in $R^d$ if and only if it satisfies one and hence both of the following equivalent conditions:

\medskip 
(i)\ it decomposes into a disjoint union of two graphs on the full set of $n$ vertices, one with
$dn-d$ edges and sparsity type, the other with ${d\choose 2}n+{d+1\choose 2}$ edges and sparsity type;

\medskip
(ii)\  it contains  the disjoint union of two graphs on the full set of $n$ vertices, one with $dn-d$ edges and sparsity type, the other with ${d\choose 2}n$ edges and sparsity type.
\end{theorem}

\medskip 
As shown in \cite{borcea:streinu:minimalPeriodic:LMS:2011}, graphs with $dn-d$ edges on $n$ vertices and with $dn-d$
sparsity type are minimally rigid graphs for the problem of rigidly connecting with
bars a set of $n$ translated bodies (i.e. parallel Cartesian frames). Thus, the presence of a
$dn-d$ sparse subgraph is related to the fact that vertex-orbits (when `frozen solid' together into a body)
provide a system of $n$ translated bodies which has to become rigid. 

\medskip 
\noindent
{\em Proof of sufficiency:} \ We begin with some combinatorial background about sparsity of type $an+b$, with $a$ and $b$ non-negative integers. Again, in this discussion, graphs may have multiple edges and loops.

\medskip 
\noindent
{\bf Loop breaking.}\ When a graph has loops, we shall describe as {\em loop breaking}
the operation of replacing a loop by an edge connecting the loop vertex with some other
vertex of the graph. It is immediate that loop breaking preserves $an+b$ sparsity. The
following lemma is the relevant converse.

\begin{lemma}\label{sparsity}
Let $G$ be a graph with $an+b$ edges on $n$ vertices, with $a$ and $b$ non-negative
integers. Suppose $G$ is $an+b$ sparse, that is, any subset of $n'$ vertices has at most
$an'+b$ edges between them. Then, there is a graph $\tilde{G}$, made of $an+b$ loops on $n$
vertices, 
which is $an+b$ sparse and yields $G$ after an adequate sequence of loop breaking.
\end{lemma}

\medskip 
\noindent
{\em Proof:}\ We claim that we can always replace an edge of $G$ between vertices $u$ and $v$
by a loop at one of these two vertices so that $an+b$ sparsity is preserved. Suppose this were
not the case. Then there are subsets of vertices $U$ and $V$, with $u$ in $U$ but not in $V$
and $v$ in $V$ but not in $U$ with $a|U|+b$ edges in $U$ and $a|V|+b$ edges in $V$.
Since the union $U\cup V$ has at most $a|U\cup V|+b$ edges, the intersection $U\cap V$ must be non-void and tight. Even so,  the union $U\cup V$ has
already $a|U\cup V|+b$ edges without counting the edge we started with. This
contradiction proves the claim and the lemma follows by iteration. Note that the resulting graph
$\tilde{G}$ must have at least $a$ loops at each vertex.\qed

\begin{corollary}
	\label{cor:anbSparse}
With the notations of the above lemma, $G$ is $an+b$ sparse if and only if it contains a
subgraph with $|V|=an$ edges on the $n$ vertices and sparsity type $an$.
\end{corollary}

\medskip 
We show now that a {\em loop breaking} operation may be performed on minimal
infinitesimally rigid frameworks. Let us assume that such a framework has a loop at vertex $j$
in the quotient graph $G/\Gamma$. Using normalized coordinates as in (\ref{eq:normalized}), the loop corresponds to an equation of the form:
\begin{equation}\label{loopRow}
 \langle \dot{\Lambda}c +A_j(q^j-q'^j), \Lambda c +q^j-q'^j\rangle =0
\end{equation}
\noindent
which involves only the difference $q^j-q'^j$. This allows the assumption $q'^j=0$ for the
origin of the bar.  We are going to replace this loop at vertex $j$ of $G/\Gamma$ by an
edge from vertex $i$ to vertex $j$. The corresponding bar has one endpoint at the origin of the body representative indexed by $i$ and the other endpoint at $kq^j$ on the body representative indexed by $j$ translated by period $\Lambda kc$, for some integer $k$. This gives an equation of the form:
\begin{equation}\label{exchangeRow}
 \langle \dot{\Lambda}kc +A_jkq^j+\dot{p}_j-\dot{p}_i, \Lambda kc +kq^j\rangle =0
\end{equation}
\noindent
which, divided by $k^2$, takes the form:
\begin{equation}\label{exchangeRowBis}
 \langle \dot{\Lambda} c +A_jq^j  +\frac{1}{k}(\dot{p}_j-\dot{p}_i), \Lambda c +q^j\rangle =0
\end{equation}
\medskip \noindent
With all other rows unchanged, the rigidity matrices for the two frameworks involve only
the exchange of row (\ref{loopRow}), where $q'_j=0$, with row (\ref{exchangeRowBis}) and for
a sufficiently large integer $k$, the rank won't be affected. 
This shows that loop breaking operations can be performed on quotients of minimally rigid graphs with preservation of minimal rigidity.

\medskip 
The sufficiency proof follows from these results.
Indeed, if $G/\Gamma$ has edges separated into a $dn-d$ tight sparse graph and a ${d\choose 2}n+{d+1\choose 2}$ tight sparse graph, we can produce minimally rigid $d$-periodic
graphs with that quotient by a two-step procedure. First, we use the $dn-d$ sparse subgraph,
but replace the remainder graph with a graph made only of loops, as guaranteed by Lemma~\ref{sparsity}. Infinitesimally rigid frameworks corresponding to this quotient structure
can be constructed in a natural way, as indicated immediately below. Then, as a second step,
we break the appropriate loops and obtain the initial graph while maintaining minimal rigidity.

\medskip 
\noindent
When the  ${d\choose 2}n+{d+1\choose 2}$ tight sparse subgraph is made entirely of loops, natural liftings of $G/\Gamma$ to minimally rigid $(G,\Gamma)$ are constructed as follows. First, we use 
${d+1\choose 2}$ loops to make rigid the lattice of periods, as in  \cite{borcea:streinu:periodicFlexibility:2010,borcea:streinu:minimalPeriodic:LMS:2011}. Note that
bars placed as periods can be interpreted as loops at any vertex of the quotient graph. We distribute them so that exactly ${d\choose 2}$ loops remain for bar assignment at each vertex class. Now in each vertex orbit we place the available ${d\choose 2}$ classes of bars so that the whole body-and-bar orbit becomes a single rigid body. Thus, with the ${d\choose 2}n+{d+1\choose 2}$ loop subgraph we have turned the system into $n$ parallel bodies and now the $dn-d$ sparse subgraph is precisely what we need to turn this system into
a single rigid body, as in \cite{borcea:streinu:minimalPeriodic:LMS:2011}.\qed

An alternative proof of sufficiency is given below after completing the necessity.

\medskip \noindent
{\em Proof of necessity:}\ We consider the $m\times u$ rigidity matrix of an infinitesimally
rigid periodic body-and-bar framework in $R^d$.  The number of rows $m={d+1\choose 2}(n-1)+d^2$ is the number of edge orbits and the number of columns $u={d+1\choose 2}n+d^2=m+{d+1\choose 2}$ is the number
of {\em `unknowns'}  in the linear constraint system for infinitesimal variations.  By assumption, the rows are linearly independent and the ${d+1\choose 2}$-dimensional kernel of our rigidity matrix consists precisely of the trivial infinitesimal rigid motions.

When we consider some subset of rows, denoted say by $F$, we have $|F|=rk(F)=u-dim(Ker(F))$.

\medskip \noindent
We denote by $n_F$  the number of vertex orbits incident to the given edge orbits and by  $\omega_F$ the number of connected components of $F$ as a set of edges in the quotient
graph. We observe that 
\begin{equation}\label{eq:kernel}
dim(Ker(F))\geq {d+1\choose 2}+(n-n_F){d+1\choose 2}+d(\omega_F-1)
\end{equation}
\noindent
where the numbers on the right hand side count the following independent infinitesimal
deformations respecting the constraint system $F$: trivial infinitesimal rigid motions, 
arbitrary infinitesimal deformations for the body representatives for non-incident vertices
and relative infinitesimal translations of the connected components in $F$. Thus, we have
a necessary {\em edge sparsity} condition:

\begin{equation}\label{eq:edgeSparsity}
|F|=u-dim(Ker(F))\leq d(n_F-\omega_F)+{d\choose 2}n_F+{d+1\choose 2}
\end{equation}

\medskip \noindent

Written in this form, the right hand side becomes very suggestive from the point of view of
{\em matroid theory} \cite{oxley:matroid:1992}. Indeed, if we consider a ground set of sufficiently many edges
and loops on $n$ vertices (so that any graph with $m$ edges can be found as a subgraph), 
then the function defined by the right hand side is an increasing submodular function with
values in the non-negative integers and determines a matroid which is, by a theorem of Pym and Perfect \cite{pym:perfect:submodular:1970}, the union of $d$ copies of the {\em graphic matroid} determined by
$f_1(F)=n_F-\omega_F$, then ${d\choose 2}$ copies of the {\em bi-cycle matroid} determined
by $f_2(F)=n_F$ and the {\em uniform matroid} determined by $f_3(F)={d+1\choose 2}$.

\noindent
A basis in a union of $d$ graphic matroids is a union of $d$ spanning trees and this is precisely 
what $dn-d$ tight sparse graphs are, while ${d\choose 2}n+{d+1\choose 2}$ tight sparse
graphs give the bases of the remaining union. Thus, the $m$ edges corresponding to
our maximal rank rigidity matrix have the claimed decomposition. \qed

This completes the proof of our main Theorem \ref{thm:minRigidity}. 

The result proven above for minimally rigid $d$-periodic body-and-bar frameworks is obviously related to the
characterization obtained in \cite{borcea:streinu:minimalPeriodic:LMS:2011} for minimally rigid $d$-periodic bar-and-joint graphs. These two types of periodic framework structures may be seen as the extreme cases $k = d$, respectively $k = 0$, of the family
obtained by articulating with bars $k$-dimensional plates, with $0\leq k\leq d$. In the full paper, we show that our proof technique can be extended to these general periodic frameworks.

\medskip \noindent
To illustrate other techniques, we present now an {\bf alternative proof of sufficiency }\ for
Theorem~\ref{thm:minRigidity}.

\medskip
Suppose the bar corresponding to $e\in E$ connects two
bodies $B_i$ and $B_j$, which are identified
with Cartesian frames $(p_i,M_i)$ and $(p_j,M_j)$ with $p_i,p_j\in R^d$ and $M_i,M_j\in SO(d)$. We normalize them (i.e., $p_i=0$ and $M_i=I_d$).
As we explained previously, a bar gives a linear equation
\begin{equation}
\label{eq:const}
\langle \dot{\Lambda}c +A_jq_j-A_iq_i+\dot{p}_j-\dot{p}_i, \Lambda c +q_j-q_i \rangle =0,
\end{equation}
where $A_i$ is a skew-symmetric matrix with $\dot{M}_i=M_iA_i$ and $q_i, q_j$ are endpoints of the bar in each frame.

Let $w_i$ be a vector in $R^{d \choose 2}$ which is obtained by aligning ${d \choose 2}$ independent entries of $A_i$.
By setting $h=q_j+\Lambda c - q_i$,
we can simply write (\ref{eq:const}) by 
\begin{equation}
\label{eq:const2}
\langle \dot{p}_j-\dot{p}_i, h \rangle+ (\langle w_j, q_j\wedge h \rangle-\langle w_i, q_i\wedge h\rangle) + \langle \dot{\Lambda} c, h\rangle =0.
\end{equation}
The rigidity matrix is then redefine as the matrix of the linear system (\ref{eq:const2}) in $\dot{p}_i, w_i$ and $\dot{\Lambda}$ over all edge representatives.
The columns of the rigidity matrix are naturally divided into three parts, the translation-part, the rotation-part, and the lattice-part.
Namely, they are the column spaces associated with $\dot{p}$, $\dot {w}$ and $\dot{\Lambda}$, respectively, and have the dimensions $dn$, ${d \choose 2}n$, and ${d+1\choose 2}$, respectively.
Also each row of the rigidity matrix is a $(dn+{d \choose 2}n+d^2)$-dimensional vector whose coordinates can be also divided into three parts.

Let $\mathbf{e}_i$ be the $i$-th element of the standard basis of $R^d$.
The following proposition gives an alternative proof of the sufficiency of our main theorem.

\begin{proposition}
\label{thm:generic_lifting}
Suppose the set of edge representatives can be partitioned into three parts:
\begin{itemize}
\item $d$ edge-disjoint spanning trees $T_1,\dots T_d$,
\item ${d \choose 2}$ edge-disjoint spanning pseudo-forests $F_{i,j}$ for $1\leq i\leq k$ and $i<j\leq d$ 
(where we may assume that in each $F_{ij}$ the edge set is oriented so that every vertex has in-degree exactly one), and 
\item $U=\{e^*_{ij}:1\leq i\leq d, 1\leq j\leq i\}$.
\end{itemize}
Define a realization by 
\begin{itemize}
\item $\Lambda=I_d$,
\item $C(e)=\mathbf{e}_i, q_j=q_i=0$ for $e\in T_i$,
\item $C(e)=\mathbf{e}_j, q_j=\mathbf{e}_i, q_i=0$ for $e\in F_{ij}$,
\item $C(e_{ij}^*)=N(\mathbf{e}_i+\mathbf{e}_j), q_i=q_j=0$ for $e_{ij}^*$ with $1\leq i\leq d$ and $i\leq j\leq d$,
\end{itemize}
where $N$ is an integer with $N>n$.
Then, the rigidity matrix takes full rank.
\end{proposition}

{\em Proof:}
Recall that $\dot{\Lambda}=[\dot{\lambda}_{ij}]$ denotes the infinitesimal motion of axes of the lattice.
Since any framework has ${d+1 \choose 2}$ trivial motions, we shall eliminate such trivial motions by restricting translations of a body by pinning down a point and 
by restricting rotations of the lattice.
Namely, let $\dot{p}_{v^*}=0$ for a vertex $v^*$ and $\dot{\lambda}_{ij}=0$ for every $(i,j)$ with $1\leq i\leq d-1$ and $1\leq j\leq i$.
To do that, we add new ${d+1 \choose 2}$ rows to the rigidity matrix in an appropriate manner. We now prove the resulting matrix has full rank.

We first note the entries of row for each $e\in F_{ij}$.
Since $q_2\wedge h=q_2\wedge (q_2+C(e)-q_1)=\mathbf{e}_j\wedge \mathbf{e}_i$  and $q_1\wedge h=0$, 
the row in the rotation-part is written as 
\[
(\Bvector{}{0},\Bvector{\cdots}{\cdots}, \Bvector{}{0},  
\Bvector{(v,i,j)}{1}, \Bvector{}{0}, \Bvector{\cdots}{\cdots}, \Bvector{}{0}) 
\] 
for each $e=uv\in F_{ij}$ (where we indexed each coordinate of a vector in $R^{{d\choose 2}n}$ by a triple $(v',i,j)$ with  $v'\in V$ and $1\leq i<j\leq d$).
Thus the set of row vectors of $\bigcup_{1\leq i<j\leq d} F_{ij}$ forms the identity matrix of size ${d \choose 2}n$ in the rotation-part of the rigidity matrix.
Moreover, since $q_1=q_2=0$ for all edges $e$ in $(\bigcup_i T_i)\cup U$, all entries of the rest of rows are zero in the rotation part. 

We thus need to show that the rest of rows vectors are linearly independent when restricting to the column space of the translation-lattice-part.
We again note the entries of each row: 
For each $e=uv\in T_i$, the row in the translation part is 
\[ 
(\Bvector{}{0},\Bvector{\cdots}{\cdots}, \Bvector{}{0},  
\Bvector{(u,i)}{-1},  \Bvector{}{0}, \Bvector{\cdots}{\cdots}, 
\Bvector{}{0},
\Bvector{(v,i)}{1}, \Bvector{}{0}, \Bvector{\cdots}{\cdots}, \Bvector{}{0}),
\]
while in the lattice part the row has the form $\mathbf{e}_i\otimes \mathbf{e}_i$, that is,
\[
(\Bvector{\dot{\lambda}_{11}}{0}, \Bvector{\cdots}{\cdots}, \Bvector{}{0}, \Bvector{\dot{\lambda}_{ii}}{1}, \Bvector{}{0}, \Bvector{\cdots}{\cdots}, \Bvector{\dot{\lambda}_{dd}}{0}).
\]
For $e_{ij}^*=uv\in U$ the corresponding row in the translation part is 
\[
(\Bvector{}{0},\Bvector{\cdots}{\cdots}, \Bvector{}{0},  
\Bvector{(u,i)}{-N},  \Bvector{}{0\ \cdots 0}, \Bvector{(u,j)}{-N}, \Bvector{}{0}, \Bvector{\cdots}{\cdots}, 
\Bvector{}{0}, 
\Bvector{(v,i)}{N}, \Bvector{}{0\ \cdots 0}, \Bvector{(v,j)}{N}, \Bvector{}{0}, \Bvector{\cdots}{\cdots}, \Bvector{}{0})
\]
while in the lattice part the row has the form $N\mathbf{e}_i\otimes N\mathbf{e}_j$, that is,
\[
(\Bvector{\dot{\lambda}_{11}}{0}, \Bvector{\cdots}{\cdots}, \Bvector{}{0}, \Bvector{\dot{\lambda}_{ii}}{N^2}, 0 \cdots \ 0, 
\Bvector{\dot{\lambda}_{ij}}{N^2}, 0 \cdots \ 0 , \Bvector{\dot{\lambda}_{ji}}{N^2}, 0 \cdots \ 0, \Bvector{\dot{\lambda}_{jj}}{N^2}\Bvector{}{0}, \Bvector{\cdots}{\cdots}, \Bvector{\dot{\lambda}_{dd}}{0}).
\]

Here is an important observation:
Suppose the endpoints of $e_{ij}^*$ are $u$ and $v$ in the quotient graph.
Then, $T_i$ has the path from $u$ to $v$, which consists of a sequence of edges $e_1,e_2,\dots, e_t$.
If we sum up the rows corresponding to these edges, the translation-lattice-part becomes
\[
(\Bvector{}{0},\Bvector{\cdots}{\cdots}, \Bvector{}{0},  
\Bvector{(u,i)}{-1},  \Bvector{}{0}, \Bvector{\cdots}{\cdots}, 
\Bvector{}{0}, 
\Bvector{(v,i)}{1}, \Bvector{}{0}, \Bvector{\cdots}{\cdots}, \Bvector{}{0}
\Bvector{\dot{\lambda}_{11}}{0}, \Bvector{\cdots}{\cdots}, \Bvector{}{0}, \Bvector{\dot{\lambda}_{ii}}{t}, \Bvector{}{0}, \Bvector{\cdots}{\cdots}, \Bvector{\dot{\lambda}_{dd}}{0}),
\]
where $t$ is the length of the path, which is less than $n$.
Even if we further multiply it by $N$, we have 
\[
(\Bvector{}{0},\Bvector{\cdots}{\cdots}, \Bvector{}{0},  
\Bvector{(u,i)}{-N},  \Bvector{}{0}, \Bvector{\cdots}{\cdots}, 
\Bvector{}{0}, 
\Bvector{(v,i)}{N}, \Bvector{}{0}, \Bvector{\cdots}{\cdots}, \Bvector{}{0}
\Bvector{\dot{\lambda}_{11}}{0}, \Bvector{\cdots}{\cdots}, \Bvector{}{0}, \Bvector{\dot{\lambda}_{ii}}{tN}, \Bvector{}{0}, \Bvector{\cdots}{\cdots}, \Bvector{\dot{\lambda}_{dd}}{0}).
\]
Since $N^2>tN$, we can safely eliminate all the non-zero entries of the translation-part in the row of $e_{ij}^*\in U$.
It is easy to observe that the submatrix induced by $T_i,\dots, T_d$ together with $d$ additional rows for eliminating trivial translations has the rank $dn$ in the translation-part 
while the submatrix induced by $e_{ij}^*$ together with ${d \choose 2}$ rows for eliminating trivial rotations has the rank $d^2$ in the lattice-part.
This shows that the rigidity matrix has full rank.


\section{Mixed plate-and-bar periodic structures}
\label{sec:mixed}

The result proven above for minimally rigid $d$-periodic body-and-bar frameworks is obviously  related to the characterization obtained in \cite{borcea:streinu:minimalPeriodic:LMS:2011} for minimally rigid $d$-periodic bar-and-joint graphs.  
Actually, these two types of periodic framework structures may be seen as the extreme cases $k=d$, respectively $k=0$, of the family obtained by articulating with bars $k$-dimensional {\em plates} in $R^d$, with $0\leq k\leq d$. A plate is understood in this context as a rigid object marked by a $k$-frame in $R^d$.

\medskip \noindent
{\bf Remark:}\ The distinction between articulations made of $(d-1)$-dimensional plates, also called {\em panels},
and those made of bodies, which are $d$-dimensional plates, is that bars have to be attached
at points of the panel. Otherwise, any $(d-1)$-frame has a unique extension to a $d$-frame with standard orientation.

\medskip \noindent
We observe that, since $k$-frames in $R^d$ can be
parametrized by $SE(d)/SO(d-k)$, the periodic setting with variable period lattice will require
for minimal rigidity the relation

$$ m=[d(k+1)-{k+1\choose 2}]n +d^2-{d+1\choose 2}=
          (dn-d)+(dk-{k+1\choose 2})n+{d+1\choose 2} $$

\noindent
and the extended version of the theorem is that a quotient graph comes from an minimally rigid $d$-periodic plate-and-bar framework (with all plates of dimension $k$) if and only if  it contains the disjoint union of two graphs on the full set of $n$ vertices, one with $dn-d$ edges and sparsity type, the other with $(dk-{k+1\choose 2})n$ edges and sparsity type. Note that, as
in Theorem~\ref{thm:minRigidity} (ii), there is no condition on the remaining ${d+1\choose 2}$ edges of $G/\Gamma$. 

\medskip
We may consider a more general setting by {\em mixing} various types of plates. 
The associated $d$-periodic graph $(G,\Gamma)$ will have vertices labeled according to the dimension of the plate they represent. This leads to a quotient graph $G/\Gamma$ with labels, or
weights $k_v\in \{0,1,\dots, d\}$ for each $v\in V/\Gamma$. With vertex orbits labeled from 1 to $n$, we write $k_i,\ i=1,...,n$ for the corresponding (dimensional) weights. Then, the same type of reasoning will establish the following general characterization of minimal rigidity:

\begin{theorem}
\label{thm:mixedMR}
A vertex weighted quotient graph $G/\Gamma$ with
  
$$ m = d(n-1)+\sum_{i=1}^n[dk_i-{k_i+1 \choose 2}]+{d+1\choose 2} $$ 

\noindent
is the quotient graph of an infinitesimally rigid mixed plate-and-bar periodic framework in $R^d$ if and only if satisfies the sparsity condition:
\begin{equation}
\label{eq:mix}
 |F| \leq d(n_F-\omega_F)+\sum_{v\in V_F}[dk_v-{k_v+1 \choose 2}]+{d+1\choose 2} 
\end{equation}
for every nonempty $F\subseteq E/\Gamma$, where $V_F$ is the set of vertices incident to $F$.
\end{theorem}

\medskip \noindent
{\bf Remark:}\ The edge sparsity formula (\ref{eq:mix}) replaces formula (\ref{eq:edgeSparsity})
of the `pure' top case $k_v=d$. Its role and combinatorial unfolding can be treated similarly.
Let $k_v'=[dk_v-{k_v+1 \choose 2}]$.
The function defined by the right hand side of (\ref{eq:mix}) determines a matroid
which is the union of $d$ copies of the graphic matroid, the uniform matroid with rank ${d+1 \choose 2}$, 
and the matroid determined by $f_4(F)=\sum_{v\in V_F} k_v'$.

\medskip
The {\em `loop breaking scenario'} used in Section~\ref{sec:mr} remains valid for the
mixed case and shows in particular that any basis of the matroid determined by $f_4$ is
obtained from the graph with $k_v'$ loops at each vertex of $V/\Gamma$ by loop breaking.

\medskip
In this sense, {\em `archetypal'} infinitesimally rigid frameworks are obtained by fixing
the periodicity lattice through ${d+1\choose 2}$ bars which are periods, then the $k_v'$
loops at the vertex orbit of $v$ are used to eliminate the relative motion of the $k_v$-frame
with respect to the lattice and turn the whole orbit of $v$ into a rigid body. Finally, the
$n$ $\Gamma$-orbits are a sysytem of translated bodies which is rendered rigid by using the 
remaining $d$ spanning trees. An alternative, more detailed argument for the proof of Theorem \ref{thm:mixedMR} is the following.

{\em Proof:} 
The necessity follows from the same argument as  Theorem~\ref{thm:minRigidity} by counting the number of independent motions.
We show that the sufficiency also follows from the loop breaking argument given in Theorem~\ref{thm:minRigidity}.

We first provide the corresponding combinatorial part. 
Let $k_v'=[dk_v-{k_v+1 \choose 2}]$ for simplicity.
By the theorem of Pym and Perfect on increasing submodular functions again, the function defined by the right hand side of (\ref{eq:mix}) determines the matroid,
which is the union of $d$ copies of the graphic matroid, the uniform matroid with rank ${d+1 \choose 2}$, 
and the matroid determined by $f_4(F)=\sum_{v\in V_F} k_v'$ which is also an increasing submodular function with values in non-negative integers.

To further decompose the last matroid, let $V_{\ell}=\{v\in V:k_v\geq \ell\}$ for $\ell=0,\dots, d$
and define $g_{\ell}(F)=|V_{\ell}\cap V_F|$.
Observe $f_4=\sum_{\ell=1}^{d-1} (d-\ell) g_{\ell}$, and hence the matroid determined by $f_4$ can be further decomposed into ${d \choose 2}$ matroids.
This gives us a decomposition of $\sum_{v\in V}k_v'$ edges into ${d\choose 2}$ subsets $F_{k\ell}$ for $1\leq \ell \leq d-1$ and $1\leq k\leq d-\ell$ such that $F_{k\ell }$ is a pseudo-forest spanning only vertices of $V_{\ell}$
(i.e., $F_{k\ell}$ is an edge set on $V_{\ell}$ with $|F_{k\ell}|=|V_{\ell}|$ such that each connected component contains exactly one cycle).
Each $F_{k\ell}$ can be oriented so that each vertex in $V_{\ell}$ has in-degree exactly one, 
and this orientation implies the following lemma, which supplies the converse direction of loop breaking operations.

\begin{lemma}
\label{lemma:converse}
Let $G$ be a graph with $\sum_{v\in V} k_v'$ edges on $n$ vertices, with $k_v'$ non-negative integers. 
Suppose $|F|\leq \sum_{v\in V_F} k_v'$ for $F\subseteq E$. 
Then, there is a graph $\tilde{G}$, made of $k_v'$ loops for each vertex $v\in V$, which yields $G$ after an adequate sequence of loop breakings.
\end{lemma}

(We should remark that Lemma~\ref{lemma:converse} actually follows from 
classical Hakimi's orientation theorem \cite{hakimi:degreesVertices:1965} of undirected graphs directly.)

We thus convert $G/\Gamma$ to a graph which consists of $dn-d$ tight sparse subgraph, ${d+1\choose 2}$ edges, and $k_v'$ loops for each $v\in V$.
As before ${d+1\choose 2}$ is further converted to ${d+1\choose 2}$ loops, and we have the base graph. 
${d+1\choose 2}$ loops are used to make rigid the lattice of periods.
In each vertex orbit $v$, we place $k_v'$ classes of bars, corresponding to $k_v'$ loops attached to $v$,  so that the whole $k_v$-dimensional plate-and-bar orbit becomes a single body.   
Thus, using  the $\sum_{v\in V}k_v'+{d+1\choose 2}$ loops we have turned the system into $n$ parallel bodies,
and the remaining $dn-d$ sparse subgraph makes this system into a single rigid body \cite{borcea:streinu:minimalPeriodic:LMS:2011}.

The proof is complete by performing a sequence of loop breaking operations we have shown in Theorem~\ref{thm:minRigidity}.
From above construction, we may assume that each $k_v$-dimensional plate in the realization is the one determined by the first $k_v$ column vectors of $I_d$.
Then the rigidity matrix is written in the same form:  
\begin{equation}\label{mixed_matrix}
 \langle \dot{\Lambda}c +A_jq_j-A_iq_i+\dot{p}_j-\dot{p}_i, \Lambda c +q_j-q_i\rangle =0
\end{equation}
over all edge representatives,
where $A_i$ is a skew-symmetric matrix with $k_v'$ non-zero entries. 
Thus we can apply the exactly same argument as the proof of Theorem~\ref{thm:minRigidity} to show that 
a loop breaking operation preserves infinitesimal rigidity.

%
\section{On minimal rigidity in terms of $(G,\Gamma)$}
\label{sec:moreMR}

The case $n=1$ treated in Section \ref{sec:mr} illustrates already the contrast between an answer given
in terms of $G/\Gamma$, which allows a generic lifting of edges to a $d$-periodic covering
graph, and an answer given in terms of the periodic graph $(G,\Gamma)$ itself, Theorem~\ref{thm:oneBodyOrbit}. We discuss here several aspects which must be implicated when
minimal rigidity is addressed directly on $(G,\Gamma)$. We assume $G$ connected and
$m={d+1\choose 2}(n-1)+d^2$.

\medskip
We reviewed in Section~\ref{sec:mr} the fact that the covering $G\rightarrow G/\Gamma$ with
covering group $\Gamma$ can be equivalently described through $G/\Gamma$ and a
surjective group homomorphism $H_1(G/\Gamma) \rightarrow \Gamma $. For a given
subset $F$ of edges in $G/\Gamma$, we shall consider the following numbers or
invariants defined in terms of the quotient graph {\em and} the surjective homomorphism
describing its covering $(G,\Gamma)$:\ $n_F,\omega_F, d_F, d_{F'}$, where $F'$ denotes
a connected component of $F$ when considered as a topological subspace of $G/\Gamma$.

\medskip
$n_F$ is the number of vertices in $G/\Gamma$ which are incident to edges in $F$;

\medskip
$\omega_F$ is the number of connected components of $F$ as a topological subspace
$F\subset G/\Gamma$;

\medskip \noindent
For each connected component $F'\subset F\subset G/\Gamma$ we consider the induced
maps

\begin{equation}
\label{eq:homology}
H_1(F') \rightarrow H_1(F) \rightarrow H_1(G/\Gamma) \rightarrow \Gamma
\end{equation}

\noindent
and define $d_F$ as the rank of the image of $H_1(F)$ in $\Gamma$, while $d_{F'}$ will
be the rank of the image of $H_1(F')$ in $\Gamma$. Obviously $d_{F'}\leq d_F$.

\medskip
Now, we may formulate a more refined estimate than (\ref{eq:kernel}) for the dimension of
the kernel of the matrix with row constraints determined by $F$:

$$ dim(Ker(F))\geq {d+1\choose 2}+(n-n_F){d+1\choose 2} + d(\omega_F-1)+ $$

\begin{equation}
\label{eq:kernelFull}
                       + \sum_{F'}{d-d_{F'}\choose 2} + [{d+1\choose 2}-{d_F+1\choose 2}]
\end{equation}

\noindent
where the new terms on the right hand side take account of the fact that when liftings
over a connected component $F'$ involve only a subspace of the lattice of periods, 
there are independent infinitesimal rotations in the orthogonal complement which
are not constrained by $F'$, while the last term refers to the unconstrained part of 
the lattice of periods.

\medskip \noindent
With $|F|$ denoting the number of edges in $F$, we have:

\begin{equation}
\label{eq:edgeSparsityFull}
|F|\leq d(n_F-\omega_F)+{d\choose 2}n_F-\sum_{F'}{d-d_{F'}\choose 2} +{d_F+1\choose 2}
\end{equation}

\noindent
The question whether this condition is also sufficient for a combinatorial characterization of
minimal rigidity is left as an {\em open problem}. We note that for $n=1$, condition (\ref{eq:edgeSparsityFull}) becomes (\ref{eq:minimal}) since $F$ is connected and 
$d_F=dim(span C(F))$.


\section{Algorithms}
\label{sec:algorithms}



We present now some algorithmic consequences of the main theoretical result, Theorem \ref{thm:minRigidity}. 
Since our Maxwell-Laman characterization decomposes the matroids associated with
minimal rigidity of periodic frameworks into smaller pieces of
well-known (graphic and bi-cycle) matroids, we can simply apply a matroid partitioning algorithm (e.g. \cite{edmonds:matroidPartition:1965}) for checking
whether a given quotient graph $G/\Gamma$ is realizable as a minimally rigid periodic framework.

We can however give a better algorithm, using condition (ii) Theorem \ref{thm:minRigidity}. It follows from a theorem of \cite{haas:lee:streinu:theran:maps:2007} that
a graph is the disjoint union of a $dn-d$ tight subgraph and a ${d
\choose 2}n$ tight subgraph  if and only if it is a ${d+1\choose
2}n-d$ tight graph. This type of sparsity is {\em  matroidal}
(i.e., the family of tight subgraphs are the bases of a matroid on the edges of the
complete graph), allowing for the use of the pebble game algorithms of \cite{streinu:lee:pebbleGames:2008}.
Thus condition (ii) can be checked in $O(n^2)$ time.
The same algorithm computes the degrees of
freedom of a generic realization of a given quotient graph, where $m$
may not be exactly ${d+1\choose 2}(n-1)+d^2$. This is done by first computing
a {\em maximum} ${d+1\choose 2}n-d$-sparse subgraph and then adding at most
${d+1 \choose 2}$ unused edges. The number of edges of the resulting
graph gives the maximum size of independent equations, its difference
from ${d+1\choose 2}(n-1)+d^2$ being the degree of freedom.

%



\subsection*{Conclusions}
\label{sec:conclusions}

Motivated by computational studies of crystalline materials, we introduced periodic body-and-bar frameworks. We gave a combinatorial characterization of those which are generically minimally rigid, obtaining the first polynomial time algorithms for recognizing them, and for computing flexibility parameters (degrees of freedom) of those which are not rigid. In the full version of the paper, attached as an appendix, we have also characterized mixed plate-and-bar periodic structures, and included additional comments on minimal rigidity in terms of $(G,\Gamma)$.

Remaining open questions include: (a) characterizing body-and-hinge and (b) body-and-pin periodic frameworks, and (c) verifying whether the {\em Molecular conjecture} \cite{katoh:tanigawa:proofMolecularConjecture:DCG:2011} holds  as well in the periodic case. Besides intrinsic theoretical interest, such results would have immediate applications, as these structures appear in modeling families of natural crystals.


\bibliographystyle{plain}

\end{document}